\documentclass[reqno,english]{amsart}
\usepackage{amsfonts,amsmath,latexsym,verbatim,amscd,mathrsfs,color,array}

\usepackage{amsmath,amssymb,amsthm,amsfonts,graphicx,color}
\usepackage{amssymb}
\usepackage{pdfsync}
\usepackage{epstopdf}
\usepackage[colorlinks=true]{hyperref}
\usepackage{subfigure}

\usepackage[utf8]{inputenc}

\makeatletter
\def\@currentlabel{2.1}\label{e:dispaa}
\def\@currentlabel{2.21}\label{e:dispau}
\def\@currentlabel{2.22}\label{e:dispav}
\def\@currentlabel{2.23}\label{e:dispaw}
\def\@currentlabel{2.24}\label{e:dispax}
\def\theequation{\thesection.\@arabic\c@equation}
\makeatother
\let\oldbibliography\thebibliography
\renewcommand{\thebibliography}[1]{%
\oldbibliography{#1}%
\setlength{\itemsep}{0pt}%
}

\renewcommand{\theequation}{\thesection.\arabic{equation}}
\newtheorem{lemma}{Lemma}[section]
\newtheorem{definition}{Definition}[section]
\newtheorem{proposition}{Proposition}[section]
\newtheorem{corollary}{Corollary}[section]
\newtheorem{remark}{Remark}[section]
\newcommand{\bremark}{\begin{remark} \em}
\newcommand{\eremark}{\end{remark} }

\newtheorem{theorem}{Theorem}[section]
\newtheorem{open problem}{Open Problem}[section]
\newtheorem{open question}{Open Question}[section]
\newtheorem{open questions on hexagonal lattice}{Open Questions on hexagonal lattice}[section]

\newcommand{\R}{{\mathbb R}}

\newcommand{\BE}{\begin{equation}}
\newcommand{\BEN}{\begin{equation*}}
\newcommand{\EE}{\end{equation}}
\newcommand{\EEN}{\end{equation*}}
\newcommand{\BL}{\begin{lemma}}
\newcommand{\EL}{\end{lemma}}
\newcommand{\BT}{\begin{theorem}}
\newcommand{\ET}{\end{theorem}}
\newcommand{\BP}{\begin{proposition}}
\newcommand{\EP}{\end{proposition}}
\newcommand{\BC}{\begin{corollary}}
\newcommand{\EC}{\end{corollary}}

\renewcommand{\Re}{\operatorname{Re}}
\renewcommand{\Im}{\operatorname{Im}}

\numberwithin{equation}{section}

\begin{document}


\title[Lattice Hexagonal Crystallization]{\bf On lattice hexagonal crystallization for non-monotone potentials}

\author{Senping Luo}

\author{Juncheng Wei}

\address[S.~Luo]{School of Mathematics and statistics, Jiangxi Normal University, Nanchang, 330022, China}
\address[J.~Wei]{Department of Mathematics, University of British Columbia, Vancouver, B.C., Canada, V6T 1Z2}

\email[S.~Luo]{luosp1989@163.com}

\email[J.~Wei]{jcwei@math.ubc.ca}

\begin{abstract}
Let $L =\sqrt{\frac{1}{\Im(z)}}\Big({\mathbb Z}\oplus z{\mathbb Z}\Big)$   where $z \in \mathbb{H}=\{z= x+ i y\;\hbox{or}\;(x,y)\in\mathbb{C}: y>0\}$ be the two dimensional lattices with unit density. Assuming that $\alpha\geq1$, we prove that
\begin{equation}\aligned\nonumber
\min_{L}\sum_{\mathbb{P}\in L, |L|=1}|\mathbb{P}|^2 e^{- \pi\alpha|\mathbb{P}|^2}
\endaligned\end{equation}
is achieved at hexagonal lattice. More generally we prove that for $\alpha \geq 1$
\begin{equation}\aligned\nonumber
\min_{L}\sum_{\mathbb{P}\in L, |L|=1}(|\mathbb{P}|^2-\frac{b}{\alpha}) e^{- \pi\alpha|\mathbb{P}|^2}
\endaligned\end{equation}
is achieved at hexagonal lattice for $b\leq\frac{1}{2\pi}$ and does not exist for $b>\frac{1}{2\pi}$.
As a consequence, we provide two classes of non-monotone potentials which lead to hexagonal crystallization among lattices.
Our results partially answer some questions raised in  \cite{Oreport, Bet2016, Bet2018, Bet2019AMP} and extend the main results in
\cite{LW2022} on minima of difference of two theta functions.
\end{abstract}

\maketitle

\setcounter{equation}{0}

\section{Introduction and main results}
\setcounter{equation}{0}

Let $L$ be a two dimensional lattice, i.e., of the form $\Big({\mathbb Z}\vec{u}\oplus {\mathbb Z}\vec{v}\Big)$, where $\vec{u}$ and $\vec{v}$ are two independent two-dimensional vectors. Many physical, chemical and number theoritic problems can be formulated to the following minimization problem on lattices:
 \begin{equation}\aligned\label{EFL}
\min_L E_f(L ), \;\;\hbox{where}\;\;E_f(L):=\sum_{\mathbb{P}\in L \backslash\{0\}} f(|\mathbb{P}|^2),\; |\cdot|\;\hbox{is the Euclidean norm on}\;\mathbb{R}^2.
\endaligned\end{equation}
See e.g. \cite{Abr,Bet2015,Bet2016,Bet2017,Bet2018a,Bet2018,Bet2019a,Bet2019,
Bet2020a,Bet2019AMP,Betermin2021JPA,Betermin2021AHP,Betermin2021a,
Betermin2021b,Betermin2022a,Bla2015,Bla2022,Che2007,Cohen,Sar2006,Radin1987,LWLJ,LW2022,
Radin1991,Serfaty2010,Serfaty2012,Serfaty2014,Serfaty2018,Sch2011,Sigal2018,Yau2019}.
The summation ranges over all the lattice points except for the origin $0$ and the function $f$ denotes the potential of the system.
 The function
$E_f(L)$ denotes the limit energy per particle of the system under the background potential $f$ over a periodical lattice $L$.

Let $ z\in \mathbb{H}:=\{z= x+ i y\;\hbox{or}\;(x,y)\in\mathbb{C}: y>0\}$. For  lattice $L$ with  unit cell area we can  use the parametrization $L =\sqrt{\frac{1}{\Im(z)}}\Big({\mathbb Z}\oplus z{\mathbb Z}\Big)$   where $z \in \mathbb{H}$. The hexagonal lattice is the lattice spanned by the two basis vectors $(1,0)$ and $(\frac{1}{2},\frac{\sqrt3}{2})$, i.e., it can be expressed by $\sqrt{\frac{A}{\frac{\sqrt3}{2}}}[\mathbb{Z}(1,0)\oplus\mathbb{Z}(\frac{1}{2},\frac{\sqrt3}{2})]$, or simply $\sqrt{\frac{A}{\frac{\sqrt3}{2}}}[\mathbb{Z}\oplus\mathbb{Z}(\frac{1}{2},\frac{\sqrt3}{2})]$, where $A$ is density/volume of the lattice.
When $A=1$, using the notation of \cite{Bet2016}, one denotes that
 \begin{equation}\aligned\label{L555}
\Lambda_1:=\sqrt{\frac{1}{\frac{\sqrt3}{2}}}[\mathbb{Z}(1,0)\oplus\mathbb{Z}(\frac{1}{2},\frac{\sqrt3}{2})]
\hbox{is the hexagonal lattice with unit density. }
\endaligned\end{equation}

 Our paper is motivated by several open questions in the minimization problem (\ref{EFL}).
In Oberwolfach report \cite{Oreport} and \cite{Bet2018}(page 3974),
B\'etermin formulated a very fundamental question, i.e,

\begin{open question}[\cite{Oreport,Bet2018}]\label{Open1}
{For an absolutely summable interaction potential $f$, what is the minimizer of $E_f(L)$ among lattices $L$$?$}

\end{open question}
Open question \ref{Open1} was further explained from physical and mathematical sides in \cite{Betermin2021JPA} as follows
\begin{itemize}
  \item {\em {\bf Physical Problem.} Assuming that strong forces (like in metals) obliges the bonds to have a
certain fixed length as well as an exact (or minimum) coordination number for each atom, what
is the crystal lattice structure with the lowest potential energy?}
  \item {\em {\bf Mathematical Problem.} For any fixed $\lambda$, what is the
minimizer of $E_f(L)$ as in $\mathcal{L}_d(\lambda)$(i.e., among lattices)?}
\end{itemize}

Open question \ref{Open1} brings about the fundamental  project on lattice minimization and crystallization, as further noted by \cite{Bet2018}, {\em "any optimality result for
$L\rightarrow E_f(L)$ supports the associated crystallization conjecture for particles through $f$."} It has plentiful applications in
solid state and statistical physics (e.g., \cite{Radin1987,Bla2015,Bla2022,Serfaty2018}), as well as analytical number theory(e.g.,
\cite{Cohen,Sar2006,Radin1991}).
In particular, hexagonal crystallization, i.e., when hexagonal(triangular) lattice is the minimizer of $E_f(L)$ among lattices $L$, attracts much more attention. The following three open questions were posed
by \cite{Bet2019AMP} on hexagonal(triangular) lattice.

\begin{open questions on hexagonal lattice}[B\'etermin-Petrache, \cite{Bet2019AMP}]\label{Open2}
\begin{itemize}
  \item  If $f(r^2)$ is not a positive superposition of Gaussians, can the triangular(hexagonal)
lattice still be a minimizer of $E_f(L)$ among lattices at any fixed density?
  \item "How negative" can the inverse Laplace transform be, while preserving
the property that the minimum $\min_{L} E_{f}(\lambda L)$ is achieved at all
$\lambda>0$ by the triangular lattice?
  \item What is the largest class of functions $f$ such that for any $\lambda>0$, the triangular lattice
is the unique minimizer of $L\rightarrow E_{f}(\lambda L)$?
\end{itemize}

\end{open questions on hexagonal lattice}

In fact, Open questions on hexagonal lattice \ref{Open2} asked that under which potentials the system admits hexagonal crystallization among lattices as proposed by B\'etermin \cite{Bet2018}.
It motivates the following more general basic question:
\begin{open question}\label{Open5}
Finding the non-monotone potentials$(f\in\mathcal{F}_d\setminus\mathcal{F}^{cm}_d)$, such that the minimum  of $L\rightarrow E_{f}( L)$ is achieved at hexagonal lattice?
\end{open question}

\begin{remark}
Here  we follow the notations in \cite{Betermin2021AHP}
\begin{equation}\aligned
\mathcal{F}_d:=\{
f: \R_+\rightarrow \R, f(r)=\int_0^\infty e^{-r t}d\mu_f(t), \mu_f(t)\in\mathcal{M}_d, |f(r)|=O(r^{-p_f})\;\hbox{as}\;
r\rightarrow\infty, p_f>\frac{1}{2}
\},
\endaligned\end{equation}
 $\mathcal{M}_d$ is the space of Radon measures on $\R_+$, and
\begin{equation}\aligned
\mathcal{F}^{cm}_d:=\{
f\in\mathcal{F}^{cm}_d, \mu_f\geq0
\}
\endaligned\end{equation}
denoting the class of admissible and completely monotone potentials.
Note that
\begin{equation}\aligned\label{CM}
\mathcal{F}^{cm}_d\;\;\hbox{is generated by Gaussian potential:}\;\; f(|\cdot|^2)=e^{-\alpha|\cdot|^2 }
\;\;\hbox{via Laplace transform}.
\endaligned\end{equation}
\end{remark}

Open question \ref{Open5} indeed asked when the hexagonal crystallization among lattices happens for non-monotone potentials.
 As
remarked by \cite{Bet2019AMP} and \cite{Betermin2021AHP} respectively,
{\em "the global optimality of a given lattice for $E_f(L)$ among all lattices can be proved rigorously in very few examples"
and "
Only few rigorous results are available on minimization of charged structures among lattices.
"
}
 The famous result by Montgomery  \cite{Mon1988}
proved that if $f(|\cdot|)=e^{-\alpha|\cdot|^2 }$ for all $\alpha>0$, then minimizer of $L\rightarrow E_{f}( L)$ is achieved at hexagonal lattice. Using the Laplace transform and Hausdorff-Bernstein-Widder representation Theorem(\cite{Bet2016}),
this leads to for the admissible completely monotone class denoted by $\mathcal{F}^{cm}_d$, the minimizer of $L\rightarrow E_{f}( L)$ is hexagonal lattice, then recover the Riesz potentials case by \cite{Ran1953,Cas1959,Dia1964,Enn1964a, Enn1964b}. There are some non-monotone potentials proved to admit hexagonal minimizer(\cite{Bet2016,Bet2018,Betermin2021a}), by Montgomery's Theorem \cite{Mon1988}) and B\'etermin's method(\cite{Bet2016}, see Theorem 2.9 in \cite{Betermin2021AHP} for a general criterion).

We are interested in the following interesting problem initiated by B\'etermin \cite{Bet2016} on hexagonal crystallization.
In Subsection 4.3 of \cite{Bet2016}(page 3249), B\'etermin asserted that the following inequality
\begin{equation}\aligned\label{PQ}
\sum_{\mathbb{P}\in \Lambda_1, |\Lambda_1|=1}|\mathbb{P}|^2 e^{- \pi\alpha|\mathbb{P}|^2}\leq \sum_{\mathbb{P}\in L, |L|=1}|\mathbb{P}|^2 e^{- \pi\alpha|\mathbb{P}|^2}
\endaligned\end{equation}
cannot be true for all $\alpha>0$, here $\Lambda_1$ is the hexagonal lattice with unit density\eqref{L555}.
This gives the following problem

\begin{open question}\label{Open6} Is there any $\alpha$ such that inequality \eqref{PQ} is true?
\end{open question}

In this paper, we provide several classes of non-monotone potentials to hexagonal crystallization among lattices.
Thereby, we give positive and partial answers to Open questions \ref{Open1}-\ref{Open6}.

We state our main results in Theorems \ref{Th1} and \ref{Th3}.
\begin{theorem}\label{Th1} Assume that $\alpha\geq1$. Consider the  following lattice minimization problem
\begin{equation}\aligned
\min_{L}\sum_{\mathbb{P}\in L, |L|=1}(|\mathbb{P}|^2-\frac{b}{\alpha}) e^{- \pi\alpha|\mathbb{P}|^2}.
\endaligned\end{equation}
There exists $b_c=\frac{1}{2\pi}$$($independent of $\alpha$$)$ such that
\begin{itemize}
  \item if $b\leq b_c$, the minimizer of the lattice energy functional is $e^{i\frac{\pi}{3}}$, which corresponds to the hexagonal lattice;

  \item if $b>b_c$, the minimizer of the lattice energy functional does not exists.
\end{itemize}

\end{theorem}
\begin{remark} In Theorems \ref{Th1} and \ref{Th3}, we provide the following two basic non-monotone potentials which lead to hexagonal crystallization among lattices.
\begin{equation}\aligned\label{Hexagon}
&(1): \hbox{Differences of Gaussian potentials},\; f(r^2)=e^{-\pi\alpha r^2}-be^{-\pi a\alpha r^2}, \;\;\alpha\geq1, a>1, b\leq\sqrt a,\;\;\\
&\;\;\;\;\;\;\;\;\;\hbox{the minimizer is hexagonal lattice};\\
&(2): \hbox{Product of Polynomial and Gaussian},\; f(r^2)=(r^2-\frac{b}{\alpha})e^{-\pi\alpha r^2}, \;\;\alpha\geq1, b\leq\frac{1}{2\pi},\\
&\;\;\;\;\;\;\;\;\;\hbox{the minimizer is hexagonal lattice}.
\endaligned\end{equation}
Therefore, we give positive answers to Open questions \ref{Open1}-\ref{Open5}.
\end{remark}

\begin{remark}

Note that the Yukawa gas on the torus$($\cite{Yau2019}, formula $(2.4)$, page 10$)$ admits the form in \eqref{EFL} with $f$ replaced by Yukawa potential. Here we consider \eqref{EFL} with $f$ replaced by two classes of non-monotone potentials \eqref{Hexagon} and ask for the shape of the torus to minimize the energy per particle \eqref{EFL}, and it turns out that the hexagonal shape of the torus wins.

\end{remark}

\begin{remark} Theorem \ref{Th1} provides a pattern for hexagonal crystallization among lattices, i.e., either admits hexagonal shape or does not exist. This is contrast to the single hexagonal shape \cite{Mon1988,Bet2016} or
the rectangular-square-rhombic-hexagonal phase transitions \cite{Bet2019,Luo2019,Luo2022,LWLJ}.

\end{remark}

A direction application of Theorem \ref{Th1}($b=0$) is the following corollary which gives partial answer to Open Question \ref{Open6}.

\begin{corollary}\label{Th2} For any $\alpha\geq1$
\begin{equation}\aligned\nonumber
\min_{L}\sum_{\mathbb{P}\in L, |L|=1}|\mathbb{P}|^2 e^{- \pi\alpha|\mathbb{P}|^2}\;\;\hbox{is achieved at hexagonal lattice}.\;\;
\endaligned\end{equation}

\end{corollary}

\begin{remark} The potential corresponds to the functional in Corollary \ref{Th2} is
$$f(r^2)=r^2 e^{-\pi\alpha r^2},$$
which is also mentioned in \cite{Can2015}$($page 1214$)$.
\end{remark}

\begin{remark} The lattice sum $\sum_{\mathbb{P}\in L, |L|=1}|\mathbb{P}|^2 e^{- \pi\alpha|\mathbb{P}|^2}$ also appeared in \cite{Cohn2022}$($Section 6, formula 6.1$)$, where the generating function is
\begin{equation}\aligned\nonumber
\tilde{F}(\tau,0)=-\frac{2\pi i\tau}{d}\sum_{x\in\Lambda_d}|x|^2 e^{\pi i |x|^2\tau}.
\endaligned\end{equation}
If $\tau=i\alpha$ be an imaginary number, then
\begin{equation}\aligned\label{Fh}
\tilde{F}(i\alpha,0)=\frac{2\pi \alpha}{d}\sum_{x\in\Lambda_d}|x|^2 e^{-\pi \alpha |x|^2}.
\endaligned\end{equation}
In dimension $d=2$, \eqref{Fh} is the functional minimized in Corollary \ref{Th2}$($up to the coefficient $\frac{2\pi \alpha}{d})$.

\end{remark}

Let
 \begin{equation}
 \label{thetas}
\theta (\alpha; z):=\sum_{\mathbb{P}\in L} e^{- \pi\alpha |\mathbb{P}|^2}=\sum_{(m,n)\in\mathbb{Z}^2} e^{-\pi\alpha \frac{\pi }{y }|mz+n|^2}
\end{equation}
be the theta function, see e.g. \cite{Mon1988,Bet2016,LW2022}.

Theorem \ref{Th1} has unexpected consequence, giving the general extension of our previous result \cite{LW2022}.
\begin{theorem}\label{Th3} Assume that $\alpha\geq1$ and $a>1$. Consider the minimization problem
\begin{equation}\aligned\nonumber
\min_{z\in\mathbb{H}}\Big(\theta(\alpha;z)-b\theta(a\alpha;z)\Big)
\endaligned\end{equation}
there exists a critical value $b_T:=\sqrt{a}$ independent of $\alpha$ such that
\begin{itemize}
  \item if $b\leq b_T$, the minimizer is $e^{i\frac{\pi}{3}}$, corresponds to hexagonal lattice;
  \item if $b>b_T$, the minimizer does not exist.
\end{itemize}
\end{theorem}
\begin{remark} Note that $\min_{z\in\mathbb{H}}\Big(\theta(\alpha;z)-b\theta(2\alpha;z)\Big), \alpha\geq1$ was solved previously
, i.e., the case $a=2$ in Theorem \ref{Th3} was proved by \cite{LW2022}, here we extend it to general $a>1$ using a different strategy.
Therefore, we provide more examples to negatively answer a Conjecture by B\'etermin $($\cite{Bet2018}, last page$)$ by Corollary \ref{Coro2}.

\end{remark}

Since $\theta(\alpha;z)=\frac{1}{\alpha}\theta(\frac{1}{\alpha};z)$(see e.g. \cite{Luo2022,Mon1988}), Theorem \ref{Th3} implies that
\begin{corollary}\label{Coro2} Assume that $\lambda\in(0,1]$ and $\beta\in(0,1)$. Then
\begin{equation}\aligned
\min_{z\in\mathbb{H}}\Big(\theta(\lambda;z)-b\theta(\beta \lambda;z)  \Big)=
\begin{cases}
\hbox{is achieved at}\;\; e^{i\frac{\pi}{3}}, &\hbox{if}\;\; b\leq\sqrt\beta,\\
\hbox{does not exist},\;\;&\hbox{if}\;\; b>\sqrt\beta.
\end{cases}
\endaligned\end{equation}

\end{corollary}

\begin{remark} In B\'etermin-Petrache \cite{Bet2019AMP}$($below Remark 1.14$)$, they noted that
{ " if we try to fix a scale constraint while minimizing $E_f(\lambda L)$ for one-well potentials, then
in general we will find different minimizers at different scales. As $\lambda$ grows, it is expected
that the minimizer changes from a triangular lattice to a rhombic one, then to a square one, then to a rectangular
and then to a degenerate rectangular one
"}. This is true for many situations$($e.g., \cite{Bet2016,Bet2019,Luo2019,Luo2022,LWLJ}$)$, however Corollary \ref{Coro2} provides a class of one-well potentials such that
the minimizer either is hexagonal one or does not exist, namely only one type minimizer(hexagonal one) at different scales.

\end{remark}

 Theorems \ref{Th1}-\ref{Th3} can be generalized by the Laplace transform (inspired by B\'etermin \cite{Bet2016}).

\begin{theorem}\label{Th4}
Let the area of two dimensional lattice $L$ be normalized to 1.
Consider the minimization problem \eqref{EFL} with potential  $f_{\alpha,P}, g_{\alpha,P}$ of the following form
\begin{equation}\aligned\label{FFGG}
f_{\alpha,P}(r):&=\int_1^\infty \Big(\big(e^{-\pi\alpha x\cdot r}-b e^{-\pi a\alpha x\cdot r}\big)\cdot P(x)\Big)dx, \;\; \alpha\geq1, b\leq \sqrt a\\
g_{\alpha,P}(r):&=\int_1^\infty \Big(\big(r\cdot x-\frac{b}{\alpha}\big)e^{-\pi\alpha x\cdot r}\cdot P(x)\Big)dx, \;\; \alpha\geq1, b\leq\frac{1}{2\pi},
\endaligned\end{equation}
where the $P(x)$ is any real function$($not necessarily continuous$)$ such that $f_{\alpha,P}(r), g_{\alpha,P}(r)$ are finite and
$$P(x)\geq0.$$

Then there exists $\beta_c=\sqrt a, \beta_s=\frac{1}{2\pi}$ independent of $P(x)$ such that
\begin{itemize}
  \item if $b\leq\beta_c$, the minimizer of $E_{f_{\alpha,P}}(L)$ exists and is always hexagonal lattice.
  \item if $b\leq\beta_s$, the minimizer of $E_{g_{\alpha,P}}(L)$ exists and is always hexagonal lattice.
\end{itemize}
\end{theorem}

A particular  application of  Theorem \ref{Th4} is the classical  Yukawa potential case.
\begin{corollary}[{\bf Yukawa potential$\{\cong P(x)\equiv1\}$ of Theorem \ref{Th4}}]\label{Coro1}
Let the area of two dimensional lattice $L$ be normalized to 1.
Consider the minimization problem \eqref{EFL} with potential
 $h_{\alpha}$ of the following form
\begin{equation}\aligned\nonumber
h_{\alpha}(r):&=\frac{e^{-\pi\alpha r}}{r}-b\frac{e^{-\pi a \alpha r}}{r}, \;\; \alpha\geq1, a>1.
\endaligned\end{equation}
Then there exists $b_{c_0}=\frac{1}{\sqrt a}$ independent of parameter $\alpha$ such that
\begin{itemize}
  \item if $b\leq b_{c_0}$, the minimizer of $E_{h_{\alpha}}(L)$ exists and is always hexagonal lattice.
\end{itemize}
\end{corollary}

\begin{remark} The first rigorous result on differences of Yukawa potential of minimizer of \eqref{EFL} is proved by
 B\'etermin \cite{Bet2016}.
  Note that we provide an effective way to prove the hexagonal crystallization among lattices under differences of Yukawa potential.
\end{remark}

Another meaningful application of Theorem \ref{Th4} is the following
\begin{corollary}[{\bf $\{\cong P(x)\equiv e^{kx}, k\leq0\}$ of Theorem \ref{Th4}}]\label{Coro1}
Let the area of two dimensional lattice $L$ be normalized to 1.
Consider the minimization problem \eqref{EFL} with potential
 $I_{\alpha}$ of the following form
\begin{equation}\aligned\nonumber
I_{\alpha}(r):&={e^{-\pi\alpha r}}\cdot\frac{r}{r+b}, \;\; \alpha\geq1, b\geq0.
\endaligned\end{equation}
Then the minimizer of $E_{I_{\alpha}}(L)$ exists and is always hexagonal lattice.
\end{corollary}

\begin{figure}
\centering
 \includegraphics[scale=0.65]{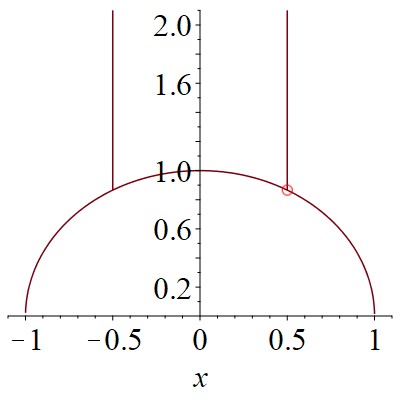}
 \caption{The hexagonal point in the fundamental domain and hexagonal shapes}
\label{f-FFF}
\end{figure}

The paper is organized as follows:
in Section 2, we provide some preliminary properties  on  the functionals and also some estimates on Jacobi theta functions.
In Section 3, we prove that the minimization on the fundamental domain can be reduced to a vertical line (see Picture \ref{f-FFF} and Theorem \ref{Th31}). In Section 4, we prove that the minimization on the vertical line can be reduced to the hexagonal point (see Picture \ref{f-FFF} and Theorem \ref{Th41}). Finally, in Section 5, we give the proof of Theorems \ref{Th1} and \ref{Th3}.


\section{Preliminaries }
\setcounter{equation}{0}

In this section we collect  some  simple symmetry properties  of the functionals and the associated fundamental domain, and also the estimates of derivatives Jacobi theta functions to be used in later sections.

Let
$
\mathbb{H}
$
 denote the upper half plane and  $\mathcal{S} $ denote the modular group
\begin{equation}\aligned\label{modular}
\mathcal{S}:=SL_2(\mathbb{Z})=\{
\left(
  \begin{array}{cc}
    a & b \\
    c & d \\
  \end{array}
\right), ad-bc=1, a, b, c, d\in\mathbb{Z}
\}.
\endaligned\end{equation}

We use the following definition of fundamental domain which is slightly different from the classical definition (see \cite{Mon1988}):
\begin{definition} [page 108, \cite{Eva1973}]
The fundamental domain associated to group $G$ is a connected domain $\mathcal{D}$ satisfies
\begin{itemize}
  \item For any $z\in\mathbb{H}$, there exists an element $\pi\in G$ such that $\pi(z)\in\overline{\mathcal{D}}$;
  \item Suppose $z_1,z_2\in\mathcal{D}$ and $\pi(z_1)=z_2$ for some $\pi\in G$, then $z_1=z_2$ and $\pi=\pm Id$.
\end{itemize}
\end{definition}

By Definition 2.1, the fundamental domain associated to modular group $\mathcal{S}$ is
\begin{equation}\aligned\label{Fd1}
\mathcal{D}_{\mathcal{S}}:=\{
z\in\mathbb{H}: |z|>1,\; -\frac{1}{2}<x<\frac{1}{2}
\}
\endaligned\end{equation}
which is open.  Note that the fundamental domain can be open. (See [page 30, \cite{Apo1976}].)

Next we introduce another group related  to the functionals $\theta(\alpha;z)$. The generators of the group are given by
\begin{equation}\aligned\label{GroupG1}
\mathcal{G}: \hbox{the group generated by} \;\;\tau\mapsto -\frac{1}{\tau},\;\; \tau\mapsto \tau+1,\;\;\tau\mapsto -\overline{\tau}.
\endaligned\end{equation}

It is easy to see that
the fundamental domain associated to group $\mathcal{G}$ denoted by $\mathcal{D}_{\mathcal{G}}$ is
\begin{equation}\aligned\label{Fd3}
\mathcal{D}_{\mathcal{G}}:=\{
z\in\mathbb{H}: |z|>1,\; 0<x<\frac{1}{2}
\}.
\endaligned\end{equation}

The following lemma characterizes the fundamental symmetries of the theta functions $\theta (s; z)$. The proof is easy so we omit it.
\begin{lemma}\label{G111} For any $s>0$, any $\gamma\in \mathcal{G}$ and $z\in\mathbb{H}$,
$\ \theta (s; \gamma(z))=\theta (s;z)$.
\end{lemma}

Let

\begin{equation}\nonumber\aligned
\mathcal{W}_b(\alpha;z):=\sum_{\mathbb{P}\in L, |L|=1}(|\mathbb{P}|^2-\frac{b}{\alpha}) e^{- \pi\alpha|\mathbb{P}|^2}.
\endaligned\end{equation}

From Lemma \ref{G111} we also have  the following invariance for $\mathcal{W}_b $.

\begin{lemma}\label{Geee}  For any $\alpha>0$ and $b\in\R$, any $\gamma\in \mathcal{G}$ and $z\in\mathbb{H}$,
$\mathcal{W}_b(\alpha; \gamma(z))=\mathcal{W}_b (\alpha;z)$.
\end{lemma}

\vskip0.05in

Next we need  some delicate analysis of the Jacobi theta function which is defined as

\begin{equation}\aligned\nonumber
\vartheta_J(z;\tau):=\sum_{n=-\infty}^\infty e^{i\pi n^2 \tau+2\pi i n z}.
 \endaligned\end{equation}
 The classical one-dimensional theta function  is given by
\begin{equation}\aligned\label{TXY}
\vartheta(X;Y):=\vartheta_J(Y;iX)=\sum_{n=-\infty}^\infty e^{-\pi n^2 X} e^{2n\pi i Y}.
 \endaligned\end{equation}
By the Poisson summation formula, it holds that
\begin{equation}\aligned\label{PXY}
\vartheta(X;Y)=X^{-\frac{1}{2}}\sum_{n=-\infty}^\infty e^{-\frac{\pi(n-Y)^2}{X}}.
 \endaligned\end{equation}

To estimate bounds of quotients of derivatives of $\vartheta(X:Y)$, we denote that

\begin{equation}\aligned\label{mmmx}
\mu(X):=\sum_{n=2}^\infty n^2 e^{-\pi(n^2-1)X},\;\;
\nu(X):=\sum_{n=2}^\infty n^4 e^{-\pi(n^2-1)X}.
\endaligned\end{equation}

We shall state a lemma which is variant of Lemmas \ref{LemmaT1} and \ref{LemmaT2} stated in the end of this section. This gives the new perspective of the estimates in Section 3.
\begin{lemma}\label{Lemma23}
\begin{itemize}
  \item  For $X>\frac{1}{5}$ and any $Y>0$, $k\in \mathbb{N}^+$
\begin{equation}\aligned\nonumber
|\frac{\vartheta_Y(X;k Y)}{\vartheta_Y(X;Y)}|\leq k\cdot\frac{1+\mu(X)}{1-\mu(X)}.
\endaligned\end{equation}
  \item For $X<\frac{\pi}{\pi+2}$ and any $Y>0$, $k\in \mathbb{N}^+$
\begin{equation}\aligned\nonumber
|\frac{\vartheta_Y(X;k Y)}{\vartheta_Y(X;Y)}|\leq k\cdot\frac{1}{\pi}e^{\frac{\pi}{4X}}.
\endaligned\end{equation}
\end{itemize}

\end{lemma}

To give the desired estimates in Section 3, we further need the following
\begin{lemma}\label{Lemma24}
\begin{itemize}
  \item For $X\geq\frac{3}{10}$ and any $Y>0$, $k\in \mathbb{N}^+$
\begin{equation}\aligned\nonumber
|\frac{\vartheta_{XY}(X;k Y)}{\vartheta_{XY}(X;Y)}|\leq k\cdot\frac{1+\nu(X)}{1-\nu(X)}.
\endaligned\end{equation}
  \item For $X\geq\frac{1}{5}$ and any $Y>0$, $k\in \mathbb{N}^+$
\begin{equation}\aligned\nonumber
|\frac{\vartheta_{XY}(X;k Y)}{\vartheta_{Y}(X;Y)}|\leq k\pi\cdot \frac{1+\nu(X)}{1-\mu(X)}.
\endaligned\end{equation}
And for $k=1$, we have the more precise bound
\begin{equation}\aligned\nonumber
|\frac{\vartheta_{XY}(X;Y)}{\vartheta_{Y}(X;Y)}|\leq\pi\cdot \frac{1+\nu(X)}{1+\mu(X)}.
\endaligned\end{equation}
\end{itemize}

\end{lemma}

\begin{proof} We first estimate $|\frac{\vartheta_{XY}(X;k Y)}{\vartheta_{XY}(X;Y)}|$ as follows
\begin{equation}\aligned\label{X1}
|\frac{\vartheta_{XY}(X;kY)}{\vartheta_{XY}(X;Y)}|
=&\frac{\sum_{n=1}^\infty n^3 e^{-\pi n^2 X}\sin(2n k\pi Y)}{\sum_{n=1}^\infty n e^{-\pi n^2 X}\sin(2n\pi Y)}\\
=&|\frac{\sin(2k\pi Y)}{\sin(2\pi Y)}|\cdot
\frac{1+\sum_{n=2}^\infty n^3 e^{-\pi (n^2-1)X}\frac{\sin(2nk\pi Y)}{\sin(2k\pi Y)}}
{1+\sum_{n=2}^\infty n^3 e^{-\pi (n^2-1)X}\frac{\sin(2n\pi Y)}{\sin(2\pi Y)}}
\endaligned\end{equation}
Then the result follows from \eqref{X1} and the following
\begin{equation}\aligned\label{X2}
|\frac{\sin(k x)}{\sin(x)}|\leq k,\;\;\hbox{for}\;\; x\in\R, k\in\mathbb{N}^+.
\endaligned\end{equation}
(The proof of (\ref{X2}) follows from a simple induction argument.)

Similar procedure applying to $\frac{\vartheta_{XY}(X;kY)}{\vartheta_{Y}(X;Y)}$ yields the desired result.

 It remains to estimate
$|\frac{\vartheta_{XY}(X;Y)}{\vartheta_{Y}(X;Y)}|$. With respect to $Y$, the function
$|\frac{\vartheta_{XY}(X;Y)}{\vartheta_{Y}(X;Y)}|$ is a periodic function with period $1$ and is symmetry about
$Y=\frac{1}{2}$. Then it suffices to consider $Y\in[0,\frac{1}{2}]$. We shall show that
\begin{equation}\aligned\label{LLL0}
\frac{\partial}{\partial Y}|\frac{\vartheta_{XY}(X;Y)}{\vartheta_{Y}(X;Y)}|\geq0\;\;\hbox{for}\;\; Y\in[0,\frac{1}{2}].
\endaligned\end{equation}
Direct computation shows that
\begin{equation}\aligned\label{LLL1}
\frac{\partial}{\partial Y}|\frac{\vartheta_{XY}(X;Y)}{\vartheta_{Y}(X;Y)}|
=\frac{D(X;Y)}{\vartheta_{Y}^2(X;Y)},
\endaligned\end{equation}
where
\begin{equation}\aligned\label{LLL2}
D(X;Y):&=
\sum_{n=1}^\infty\sum_{m=1}^\infty D_{n,m}(X;Y),\\
D_{n,m}(X;Y):&= nm(n^2-m^2)e^{-\pi(m^2+n^2-2)X}(\frac{\sin(2n\pi Y)}{\sin(2\pi Y)})'\cdot(\frac{\sin(2m\pi Y)}{\sin(2\pi Y)})
\endaligned\end{equation}
We further split the double sum into fours parts as follows
\begin{equation}\aligned\label{LLL3}
\sum_{n=1,2}\sum_{m=1,2}+\sum_{n=1,2,m\geq3}+\sum_{m=1,2,n\geq3}+\sum_{n\geq3, m\geq3}.
\endaligned\end{equation}
Note that by a direct simplification, one has
\begin{equation}\aligned\label{LLL4}
\sum_{n=1,2}\sum_{m=1,2}D_{n,m}(X;Y)=-24\pi e^{-3\pi X}\sin(2\pi Y)
\endaligned\end{equation}
and
\begin{equation}\aligned\label{LLL5}
\frac{D_{n,m}(X;Y)}{24\pi e^{-3\pi X}\sin(2\pi Y)}
=\frac{1}{24\pi}nm(n^2-m^2) e^{-\pi(m^2+n^2-5)X}\frac{1}{\sin(2\pi Y)}\big(\frac{\sin(2 n\pi Y)}{\sin(2\pi Y)}\big)'
\frac{\sin(2 m\pi Y)}{\sin(2\pi Y)}.
\endaligned\end{equation}
By Lemma \ref{Lemma27}
\begin{equation}\aligned\label{LLL6}
|\frac{D_{n,m}(X;Y)}{24\pi e^{-3\pi X}\sin(2\pi Y)}|
\leq\frac{1}{36}n^2m^2|n^2-m^2|(n^2-1) e^{-\pi(m^2+n^2-5)X}.
\endaligned\end{equation}

Therefore, by \eqref{LLL1}, \eqref{LLL2}, \eqref{LLL3}, \eqref{LLL4}, and \eqref{LLL6},

\begin{equation}\aligned\label{LLL7}
\frac{D(X;Y)}{-24\pi e^{-3\pi X}\sin(2\pi Y)}=&1-\frac{1}{36}\big(\sum_{n=1,2,m\geq3}+\sum_{m=1,2,n\geq3}+\sum_{n\geq3, m\geq3}\big)\\
&\;\;n^2m^2|n^2-m^2|(n^2-1) e^{-\pi(m^2+n^2-5)X}.
\endaligned\end{equation}
The right hand side of \eqref{LLL7} is hence is positive if $X>0.21$. Note that $1-\nu(X)>0$ if $X>0.2989938127\cdot$.

\eqref{LLL7} and \eqref{LLL1} prove \eqref{LLL0}.
By \eqref{LLL0},
\begin{equation}\aligned\nonumber
|\frac{\vartheta_{XY}(X;Y)}{\vartheta_{Y}(X;Y)}|\leq|\frac{\vartheta_{XY}(X;0)}{\vartheta_{Y}(X;0)}|,
\endaligned\end{equation}
which gives the desired result.

\end{proof}

We shall establish the following estimates which are useful in the next section.
\begin{lemma}\label{Lemma25} For $X\leq\frac{1}{2}$ and any $Y>0$, $k\in \mathbb{N}^+$
\begin{equation}\aligned\nonumber
|\frac{\vartheta_{XY}(X;Y)}{\vartheta_Y(X;Y)}|\leq \frac{3}{2}X^{-1}(1+\frac{\pi}{6}\frac{1}{X}).
\endaligned\end{equation}

\begin{equation}\aligned\nonumber
|\frac{\vartheta_{XY}(X;k Y)}{\vartheta_Y(X;Y)}|\leq \frac{3k}{2\pi}X^{-1}(1+\frac{\pi}{6}\frac{1}{X})e^{\frac{\pi}{4X}}.
\endaligned\end{equation}

\end{lemma}

\begin{proof}
By \eqref{PXY}, after a simple calculation, one has
\begin{equation}\aligned\label{PXY1}
\vartheta_{XY}(X;Y)=\pi X^{-\frac{7}{2}}
\Big(
-3X\sum_{n\in \mathbb{Z}}(n-Y) e^{-\frac{\pi(n-Y)^2}{X}}+2\pi\sum_{n\in \mathbb{Z}}(n-Y)^3 e^{-\frac{\pi(n-Y)^2}{X}}
\Big).
\endaligned\end{equation}
Then by \eqref{PXY1} and \eqref{PXY}, one has
\begin{equation}\aligned\label{PXY2}
|\frac{\vartheta_{XY}(X;Y)}{\vartheta_Y(X;Y)}|
=|\frac{3}{2}X^{-1}\cdot\Big((1-\frac{2\pi}{3X}
\cdot\frac{\sum_{n\in\mathbb{Z}}(n-Y)^3e^{-\frac{\pi(n-Y)^2}{X}}}{\sum_{n\in\mathbb{Z}}(n-Y)e^{-\frac{\pi(n-Y)^2}{X}}}\Big)|
\endaligned\end{equation}
The first part of Lemma \ref{Lemma25} follows by \eqref{PXY2} and Lemma \ref{Lemma26}.
To prove the second part of Lemma \ref{Lemma25}, one uses the following deformation
\begin{equation}\aligned\nonumber
|\frac{\vartheta_{XY}(X;kY)}{\vartheta_Y(X;Y)}|
=|\frac{\vartheta_{XY}(X;kY)}{\vartheta_Y(X;kY)}|\cdot|\frac{\vartheta_{Y}(X;kY)}{\vartheta_Y(X;Y)}|.
\endaligned\end{equation}
Therefore, the second part of Lemma \ref{Lemma25} follows by Lemma \ref{Lemma23} and the first part of Lemma \ref{Lemma25}.

\end{proof}

In Lemmas \ref{Lemma26} and \ref{Lemma27}, we provide two estimates used in Lemma \ref{Lemma25}.

\begin{lemma}\label{Lemma26}
\begin{equation}\aligned\nonumber
\sup_{X\in(0,\frac{1}{2}], Y\in\R}|\frac{\sum_{n\in\mathbb{Z}}(n-Y)^3e^{-\frac{\pi(n-Y)^2}{X}}}{\sum_{n\in\mathbb{Z}}(n-Y)e^{-\frac{\pi(n-Y)^2}{X}}}|
\leq\frac{1}{4}.
\endaligned\end{equation}

\end{lemma}

\begin{proof} Let $a=\frac{1}{X}$ and $f(a,Y):=\frac{\sum_{n\in\mathbb{Z}}(n-Y)^3e^{-a{\pi(n-Y)^2}}}{\sum_{n\in\mathbb{Z}}(n-Y)e^{-a{\pi(n-Y)^2}}}$,
then $a\geq2$. By direct checking, one has
\begin{equation}\aligned\nonumber
f(a,Y+1)=f(a,Y),\;\; f(a,1-Y)=f(a,Y).
\endaligned\end{equation}
Then it reduces to consider $f(a,Y)$ for $a\geq2$ and $Y\in[0,\frac{1}{2}]$. It suffices to prove that
$$\sup_{a\geq2, Y\in[0,\frac{1}{2}]}|f(a,Y)|\leq\frac{1}{4}.$$
After long computations we omit the details here, one has
\begin{equation}\aligned\nonumber
f_Y(a,Y)>0\;\;\hbox{for}\;\; a\geq2, Y\in[0,\frac{1}{2}].
\endaligned\end{equation}
It follows that for $a\geq2$
\begin{equation}\aligned\label{fa0}
\max_{ Y\in[0,\frac{1}{2}]}|f(a,Y)|=\max\{|f(a,Y=0)|,|f(a,Y=\frac{1}{2})|\}.
\endaligned\end{equation}
By L'Hospital's rule,
\begin{equation}\aligned\label{fa1}
|f(a,Y=0)|=\frac{\sum_{n\in\mathbb{Z}}
(2a\pi n^4-3n^2) e^{-a\pi n^2}
}{\sum_{n\in\mathbb{Z}}(2a\pi n^2-1) e^{-a\pi n^2}}
\leq\frac{4a\pi \sum_{n=1}^\infty n^4 e^{-a\pi n^2}}{1-4a\pi\sum_{n=1}^\infty n^2 e^{-a\pi n^2}}\leq0.05\;\;\hbox{for}\;\;a\geq2
\endaligned\end{equation}
and
\begin{equation}\aligned\label{fa2}
|f(a,Y=\frac{1}{2})|=&\frac{\sum_{n\in\mathbb{Z}}
(2a\pi (n-\frac{1}{2})^4-3(n-\frac{1}{2})^2) e^{-a\pi n^2}
}{\sum_{n\in\mathbb{Z}}(2a\pi (n-\frac{1}{2})^2-1) e^{-a\pi n^2}}
=\frac{1}{4}\cdot\frac{a\pi-6}{a\pi-2}\cdot\frac{1+\sigma_{a,1}}{1+\sigma_{a,2}}\\
=&
\frac{1}{4}\Big(1-(\frac{4}{a\pi-2}-\frac{\sigma_{a,1}-\sigma_{a,2}}{1+\sigma_{a,2}})
-\frac{\sigma_{a,1}-\sigma_{a,2}}{1+\sigma_{a,2}}\cdot\frac{4}{a\pi-2}\Big).
\endaligned\end{equation}
Here
\begin{equation}\aligned\nonumber
\sigma_{a,1}:&=\sum_{n=2}^\infty\frac{2a\pi (n-\frac{1}{2})^4-3(n-\frac{1}{2})^2}{\frac{a\pi}{8}-\frac{3}{4}} e^{-a\pi(n^2-n)},\\
\sigma_{a,2}:&=\sum_{n=2}^\infty\frac{2a\pi (n-\frac{1}{2})^2-1}{\frac{a\pi}{2}-1} e^{-a\pi(n^2-n)}.
\endaligned\end{equation}
Note that $\sigma_{a,1}>\sigma_{a,2}$ and
\begin{equation}\aligned\label{fa33}
a\Big(\frac{4}{a\pi-2}-\frac{\sigma_{a,1}-\sigma_{a,2}}{1+\sigma_{a,2}})\Big)\geq1>0
\endaligned\end{equation}
by a direct computation. It follows from \eqref{fa2} and \eqref{fa33} that
\begin{equation}\aligned\label{fa2a}
|f(a,Y=\frac{1}{2})|\leq\frac{1}{4}\;\;\hbox{for}\;\;a\geq2.
\endaligned\end{equation}
The bound $\frac{1}{4}$ in \eqref{fa2a} is sharp since it is approached asymptotically as $a\rightarrow\infty$ by \eqref{fa2}.
By \eqref{fa0}, \eqref{fa1} and \eqref{fa2a}, the proof is complete.

\end{proof}

The following Lemma \ref{Lemma27} is elementary and probably known in calculus, however we have not found a reference for it and so we  give the details here.

\begin{lemma}\label{Lemma27} For $n\in \mathbb{N}^+$, it holds that
\begin{equation}\aligned\nonumber
|\frac{1}{\sin(2\pi Y)}\big(\frac{\sin(2 n\pi Y)}{\sin(2\pi Y)}\big)'|\leq
C(n), \;\;\hbox{for}\;\; Y\in\R,
\endaligned\end{equation}
where $C(n)=\frac{2\pi}{3}(n-1)n(n+1)$. The upper bound $C(n)$ is sharp and is attained at $Y=k\pi,k\in N$.

\end{lemma}

\begin{proof}
Since
\begin{equation}\aligned\nonumber
\frac{\sin(2(n+1)\pi Y)}{\sin(2\pi Y)}&=\frac{\sin(2n\pi Y)\cos(2\pi Y)+\cos(2n\pi Y)\sin(2\pi Y)}{\sin(2\pi Y)}\\
&=\cos(2\pi Y)\cdot \frac{\sin(2n\pi Y)}{\sin(2\pi Y)}+\cos(2n\pi Y).
\endaligned\end{equation}
It follows that
\begin{equation}\aligned\label{L1}
\frac{1}{\sin(2\pi Y)}\big(\frac{\sin(2 (n+1)\pi Y)}{\sin(2\pi Y)}\big)'=
\cos(2\pi Y)\cdot\frac{1}{\sin(2\pi Y)}\big(\frac{\sin(2 n\pi Y)}{\sin(2\pi Y)}\big)'
-2(n+1)\pi\frac{\sin(2n\pi Y)}{\sin(2\pi Y)}.
\endaligned\end{equation}
Let
\begin{equation}\aligned\nonumber
a_n:=|\frac{1}{\sin(2\pi Y)}\big(\frac{\sin(2 n\pi Y)}{\sin(2\pi Y)}\big)'|\;\;\hbox{for short}.
\endaligned\end{equation}
Then by \eqref{L1},
\begin{equation}\aligned\label{L2}
a_{n+1}-a_{n}\leq 2\pi(n+n^2).
\endaligned\end{equation}
Here the inequality $|\frac{\sin(2n\pi Y)}{\sin(2\pi Y)}|\leq n$ is used.
By \eqref{L2},
\begin{equation}\aligned\nonumber
a_n\leq\sum_{k=1}^{n-1}(a_{k+1}-a_{k})+a_1\leq\sum_{k=1}^{n-1}2\pi(k+k^2)=\frac{2\pi}{3}(n-1)n(n+1),
\endaligned\end{equation}
which yields the desired result.

\end{proof}

The following Lemmas \ref{LemmaT1} and \ref{LemmaT2} are proved in  \cite{Luo2022}.
\begin{lemma}\label{LemmaT1}\cite{Luo2022}. Assume $X>\frac{1}{5}$. If $\sin(2\pi Y)>0$, then
\begin{equation}\aligned\nonumber
-\overline\vartheta(X)\sin(2\pi Y)\leq\frac{\partial}{\partial Y}\vartheta(X;Y)\leq-\underline\vartheta(X)\sin(2\pi Y).
 \endaligned\end{equation}
If $\sin(2\pi Y)<0$, then
\begin{equation}\aligned\nonumber
-\underline\vartheta(X)\sin(2\pi Y)\leq\frac{\partial}{\partial Y}\vartheta(X;Y)\leq-\overline\vartheta(X)\sin(2\pi Y).
 \endaligned\end{equation}
Here
\begin{equation}\aligned\nonumber
\underline\vartheta(X):=4\pi e^{-\pi X}(1-\mu(X)), \;\; \overline\vartheta(X):=4\pi e^{-\pi X}(1+\mu(X)),
 \endaligned\end{equation}
and
\begin{equation}\label{mmmx}
\mu(X):=\sum_{n=2}^\infty n^2 e^{-\pi(n^2-1)X}.
\end{equation}

\end{lemma}


\begin{lemma}\label{LemmaT2}\cite{Luo2022}.
Assume $X<\min\{\frac{\pi}{\pi+2},\frac{\pi}{4\log\pi}\}=\frac{\pi}{\pi+2}$. If $\sin(2\pi Y)>0$, then
\begin{equation}\aligned\nonumber
-\overline\vartheta(X)\sin(2\pi Y)\leq\frac{\partial}{\partial Y}\vartheta(X;Y)\leq-\underline\vartheta(X)\sin(2\pi Y).
 \endaligned\end{equation}
If $\sin(2\pi Y)<0$, then
\begin{equation}\aligned\nonumber
-\underline\vartheta(X)\sin(2\pi Y)\leq\frac{\partial}{\partial Y}\vartheta(X;Y)\leq-\overline\vartheta(X)\sin(2\pi Y).
 \endaligned\end{equation}
Here
\begin{equation}\aligned\nonumber
\underline\vartheta(X):=\pi e^{-\frac{\pi}{4X}}X^{-\frac{3}{2}};\;\; \overline\vartheta(X):=X^{-\frac{3}{2}}.
 \endaligned\end{equation}

\end{lemma}


\section{The horizontal monotonicity}

\setcounter{equation}{0}
Let $\mathcal{D}_{\mathcal{G}}:=\{
z\in\mathbb{H}: |z|>1,\; 0<x<\frac{1}{2}
\}$ be the fundamental domain associated to the group $\mathcal{G}$. Define the vertical line
\begin{equation}\aligned\label{Gfff}
\Gamma:=\{
z\in\mathbb{H}: \Re(z)=\frac{1}{2},\; \Im(z)\geq\frac{\sqrt3}{2}
\}.
\endaligned\end{equation}
See Picture \ref{f-FFF}.

Define
\begin{equation}\label{Wbbb}\aligned
\mathcal{W}_b(\alpha;z):=\sum_{\mathbb{P}\in L, |L|=1}(|\mathbb{P}|^2-\frac{b}{\alpha}) e^{- \pi\alpha|\mathbb{P}|^2}.
\endaligned\end{equation}
We use the parametrization $L =\sqrt{\frac{1}{\Im(z)}}\Big({\mathbb Z}\oplus z{\mathbb Z}\Big)$   where $z \in \mathbb{H}$, then
an explicit expression of $\mathcal{W}_b(\alpha;z)$ based on double infinite sum is
\begin{equation}\aligned
\mathcal{W}_b(\alpha;z)=\sum_{(m,n)\in\mathbb{Z}^2 }\big(
\frac{1 }{y }|mz+n|^2-\frac{b}{\alpha}
\big) e^{-\alpha \frac{\pi }{y }|mz+n|^2}.
\endaligned\end{equation}

The statement of Theorem \ref{Th1} is equivalent to
\begin{theorem}\label{Th30} Assume that $\alpha\geq1$. Then
\begin{equation}\aligned
\min_{z\in\mathbb{H}}\mathcal{W}_{b}(\alpha;z)=
\begin{cases}
\hbox{is achieved at}\;\; e^{i\frac{\pi}{3}}, &\hbox{if}\;\; b\leq\frac{1}{2\pi},\\
\hbox{does not exist},\;\;&\hbox{if}\;\; b>\frac{1}{2\pi}.
\end{cases}
\endaligned\end{equation}

\end{theorem}

We first have a comparison principle:
\begin{lemma}\label{Lemma1} Assume that $\alpha\geq1$. If
\begin{equation}\aligned
\min_{z\in\mathbb{H}}\mathcal{W}_{b_0}(\alpha;z)\;\;\hbox{is achieved at}\;\; e^{i\frac{\pi}{3}}.
\endaligned\end{equation}
Then for $b\leq b_0$,
\begin{equation}\aligned
\min_{z\in\mathbb{H}}\mathcal{W}_{b}(\alpha;z)\;\;\hbox{is stil achieved at}\;\; e^{i\frac{\pi}{3}}.
\endaligned\end{equation}

\end{lemma}

\begin{proof} For $b\leq b_0$, we use the deformation
\begin{equation}\aligned\label{Wdeform}
\mathcal{W}_{b}(\alpha;z)=\mathcal{W}_{b_0}(\alpha;z)+\frac{b_0-b}{\alpha}\cdot\theta(\alpha;z).
\endaligned\end{equation}
The result follows from \eqref{Wdeform} and the fact that
\begin{equation}\aligned
\min_{z\in\mathbb{H}}\theta(\alpha;z)\;\;\hbox{is achieved at}\;\; e^{i\frac{\pi}{3}}
\endaligned\end{equation}
by \cite{Mon1988}.

\end{proof}
To prove the first part of Theorem \ref{Th30}, by Lemma \ref{Lemma1}, we only need to prove the borderline case when $b=\frac{1}{2\pi}$
\begin{theorem}\label{Th30a} Assume that $\alpha\geq1$. Then
\begin{equation}\aligned
\min_{z\in\mathbb{H}}\mathcal{W}_{\frac{1}{2\pi}}(\alpha;z)\;\;
\hbox{is achieved at}\;\; e^{i\frac{\pi}{3}}.
\endaligned\end{equation}
\end{theorem}

The proof of Theorem \ref{Th30a} consists of two parts. In this section, we aim to prove the first part, namely

\begin{theorem}\label{Th31} Assume that $\alpha\geq1$. Then for $b=\frac{1}{2\pi}$,
\begin{equation}\aligned
\min_{z\in\mathbb{H}}\mathcal{W}_{\frac{1}{2\pi}}(\alpha;z)
=\min_{z\in\overline{\mathcal{D}_{\mathcal{G}}}}\mathcal{W}_{\frac{1}{2\pi}}(\alpha;z)
=\min_{z\in\Gamma}\mathcal{W}_{\frac{1}{2\pi}}(\alpha;z),
\endaligned\end{equation}
where $\Gamma$ is a vertical line defined at \eqref{Gfff}.
\end{theorem}

The proof of Theorem \ref{Th31} is based on the following horizontal monotonicity result:

\begin{theorem}\label{Th32} Assume that $\alpha\geq1$. Then for $b=\frac{1}{2\pi}$,
\begin{equation}\aligned
\frac{\partial}{\partial x}\mathcal{W}_{\frac{1}{2\pi}}(\alpha;z)<0,\;\;\hbox{for}\;\;z\in\mathcal{D}_{\mathcal{G}}.
\endaligned\end{equation}
\end{theorem}

In the rest of this section, we prove Theorem \ref{Th32}.

\subsection{The estimates}

We first provide en exponential expansion of $\mathcal{W}_b(\alpha;z)$, which is useful in our estimates.

\begin{lemma}\label{Lemma32} We have the following exponential expansion of $\mathcal{W}_b(\alpha;z)$:
\begin{equation}\aligned
\mathcal{W}_b(\alpha;z)=&\frac{1}{\pi}\alpha^{-\frac{5}{2}}y^{\frac{3}{2}}\cdot\Big(
\frac{1}{2}(1-2\pi b)\cdot\frac{\alpha}{y}\sum_{n\in\mathbb{Z}} e^{-\alpha\pi y n^2}\vartheta(\frac{y}{\alpha};nx)\\
&\;\;+\pi \alpha^2\cdot\sum_{n\in\mathbb{Z}} n^2e^{-\alpha\pi y n^2}\vartheta(\frac{y}{\alpha};nx)
+\sum_{n\in\mathbb{Z}} e^{-\alpha\pi y n^2}\vartheta_X(\frac{y}{\alpha};nx)
\Big).
\endaligned\end{equation}
In particular, we have the exponential expansion of $\mathcal{W}_b(\alpha;z)$ when $b$ equals to the borderline case $\frac{1}{2\pi}$
\begin{equation}\aligned\label{Wpi}
\mathcal{W}_{\frac{1}{2\pi}}(\alpha;z)=&\frac{1}{\pi}\alpha^{-\frac{5}{2}}y^{\frac{3}{2}}\cdot\Big(
\pi \alpha^2\cdot\sum_{n\in\mathbb{Z}} n^2e^{-\alpha\pi y n^2}\vartheta(\frac{y}{\alpha};nx)
+\sum_{n\in\mathbb{Z}} e^{-\alpha\pi y n^2}\vartheta_X(\frac{y}{\alpha};nx)
\Big).
\endaligned\end{equation}

\end{lemma}

The proof of Lemma \ref{Lemma32} is based on Lemmas \ref{Lemma33} and \ref{Lemma34}.

The following Lemma is based on an observation of the structure of $\mathcal{W}_{b}(\alpha;z)$\eqref{Wbbb}.
\begin{lemma}\label{Lemma33} We have the structure of $\mathcal{W}_{b}(\alpha;z)$.
\begin{equation}\aligned
\mathcal{W}_{b}(\alpha;z)=-\frac{1}{\pi}\frac{\partial}{\partial \alpha}\theta(\alpha;z)-\frac{b}{\alpha}\theta(\alpha:z).
\endaligned\end{equation}
\end{lemma}
\begin{proof} Note that
\begin{equation}\label{DDD}\aligned
\sum_{\mathbb{P}\in L, |L|=1}|\mathbb{P}|^2 e^{- \pi\alpha|\mathbb{P}|^2}=
-\frac{1}{\pi}\frac{\partial}{\partial \alpha}\sum_{\mathbb{P}\in L, |L|=1}e^{- \pi\alpha|\mathbb{P}|^2}.
\endaligned\end{equation}
The Lemma follows from \eqref{DDD} and \eqref{Wbbb}.

\end{proof}

The following Lemma is used in \cite{Luo2022,Mon1988}.
\begin{lemma}\cite{Luo2022,Mon1988}\label{Lemma34} We have the expansion of $\theta (\alpha;z)$.
\begin{equation}\aligned\nonumber
\theta (\alpha;z)
&=\sqrt{\frac{y}{\alpha}}\sum_{n\in\mathbb{Z}}e^{-\alpha \pi y n^2}\vartheta(\frac{y}{\alpha};nx)\\
&=2\sqrt{\frac{y}{\alpha}}\sum_{n=1}^\infty e^{-\alpha \pi y n^2}\vartheta(\frac{y}{\alpha};nx)+\sqrt{\frac{y}{\alpha}}\vartheta(\frac{y}{\alpha};0).
\endaligned\end{equation}

\end{lemma}

We shall also state an consequence of Lemma \ref{Lemma33}, which explains partially why we split several cases in our proof of Theorem \ref{Th32}.

\begin{lemma}\label{Lemma35} For $b=\frac{1}{2\pi}$ and any $z\in\mathbb{H}$,
\begin{equation}\aligned
\mathcal{W}_{\frac{1}{2\pi}}(1;z)=0.
\endaligned\end{equation}
\end{lemma}

\begin{proof} Since
\begin{equation}\aligned\label{Thaaa}
\theta(\frac{1}{\alpha};z)=\alpha\cdot\theta(\alpha;z)
\endaligned\end{equation}
by Fourier transform, see e.g. \cite{Luo2022,Mon1988}.
Taking derivative with respect to $\alpha$ on both sides of \eqref{Thaaa}, one gets
\begin{equation}\aligned\label{Thaaa1}
-\frac{1}{\alpha^2}\frac{\partial}{\partial \alpha}\theta(\frac{1}{\alpha};z)-\alpha\frac{\partial}{\partial \alpha}\theta(\alpha;z)
=\theta(\alpha;z).
\endaligned\end{equation}
Evaluating $\alpha=1$ at \eqref{Thaaa1}, one has
\begin{equation}\aligned\label{Thaaa2}
\Big(-\frac{\partial}{\partial \alpha}\theta(1;z)-\frac{1}{2}\theta(1;z)\Big)=0.
\endaligned\end{equation}
\eqref{Thaaa2} and Lemma \ref{Lemma33} give the result.
\end{proof}

By Lemma \ref{Lemma35}, one has
\begin{equation}\aligned\label{Eq319}
\frac{\partial}{\partial x}\mathcal{W}_{\frac{1}{2\pi}}(1;z)=0\;\;\hbox{for}\;\; z\in\mathcal{D}_{\mathcal{G}}.
\endaligned\end{equation}
Therefore, given by \eqref{Eq319}, to prove Theorem \ref{Th32}, we split the proof into two cases, namely {\bf case a: $a\in[1,1.2]$, case b: $a\in[1.2,\infty)$.}

We first give the proof of {\bf case b: $a\in[1.2,\infty)$.}
A direct consequence by taking derivative of \eqref{Wpi} in Lemma \ref{Lemma32} is
\begin{lemma}\label{Lemma36} We have the expansion of $-\frac{\partial}{\partial x}\mathcal{W}_{\frac{1}{2\pi}}(\alpha;z)$
\begin{equation}\aligned\nonumber
-\frac{\partial}{\partial x}\mathcal{W}_{\frac{1}{2\pi}}(\alpha;z)=&\frac{1}{\pi}\alpha^{-\frac{5}{2}}y^{\frac{3}{2}}\cdot\Big(
\pi \alpha^2\cdot\sum_{n\in\mathbb{Z}} n^3e^{-\alpha\pi y n^2}\cdot\big(-\vartheta_Y(\frac{y}{\alpha};nx)\big)
+\sum_{n\in\mathbb{Z}} ne^{-\alpha\pi y n^2}\big(-\vartheta_{XY}(\frac{y}{\alpha};nx)\big)
\Big).
\endaligned\end{equation}
\end{lemma}
To prove {\bf case b}, we use the deformation of $\mathcal{W}_{\frac{1}{2\pi}}(\alpha;z)$ by Lemma \ref{Lemma36}.
\begin{lemma}\label{Lemma37} We have the quotient expansion of $-\frac{\partial}{\partial x}\mathcal{W}_{\frac{1}{2\pi}}(\alpha;z)$
\begin{equation}\aligned\nonumber
-\frac{\partial}{\partial x}\mathcal{W}_{\frac{1}{2\pi}}(\alpha;z)=\frac{2}{\pi}\alpha^{-\frac{5}{2}}y^{\frac{3}{2}}
\big(-\vartheta_Y(\frac{y}{\alpha};x)\big)\cdot e^{-2\pi y}
\cdot
\Big(
\pi \alpha^2\cdot\big(1+\sum_{n=2}^\infty n^3e^{-\alpha\pi y (n^2-1)}\cdot\frac{\vartheta_Y(\frac{y}{\alpha};nx)}{\vartheta_Y(\frac{y}{\alpha};x)}\\
+ \frac{\vartheta_{XY}(\frac{y}{\alpha};x)}{\vartheta_Y(\frac{y}{\alpha};x)}
+\sum_{n=2}^\infty ne^{-\alpha\pi y (n^2-1)}\frac{\vartheta_{XY}(\frac{y}{\alpha};nx)}{\vartheta_{Y}(\frac{y}{\alpha};x)}
\Big).
\endaligned\end{equation}
Here for $\frac{y}{\alpha}>0$ and $x\in\R$,
\begin{equation}\aligned\label{Tp100}
-\vartheta_Y(\frac{y}{\alpha};x)>0.
\endaligned\end{equation}
\end{lemma}
\begin{remark} \eqref{Tp100} follows by Lemmas \ref{LemmaT1} and \ref{LemmaT2}.

\end{remark}

Based on the deformation in Lemma \ref{Lemma37} and quotient estimates of derivatives of theta function in Section 2, we are ready to prove
case b of Theorem \ref{Th32}. Namely, we are going to prove that
\begin{proposition}\label{Th32a} Assume that $\alpha\geq1.1$. Then for $b=\frac{1}{2\pi}$,
\begin{equation}\aligned
\frac{\partial}{\partial x}\mathcal{W}_{\frac{1}{2\pi}}(\alpha;z)<0,\;\;\hbox{for}\;\;z\in\mathcal{D}_{\mathcal{G}}.
\endaligned\end{equation}
\end{proposition}

\begin{proof} For convenience for stating the estimates, we denote that
\begin{equation}\aligned
\mathcal{C}(\alpha,x,y):=\frac{2}{\pi}\alpha^{-\frac{5}{2}}y^{\frac{3}{2}}
\big(-\vartheta_Y(\frac{y}{\alpha};x)\big)\cdot e^{-2\pi y}.
\endaligned\end{equation}
Here
\begin{equation}\aligned
\mathcal{C}(\pi,\alpha,y)>0
\endaligned\end{equation}
by \eqref{Tp100}.
Note that $z\in\mathcal{D}_{\mathcal{G}}$ implies that $y\geq\frac{\sqrt3}{2}$.

 We further split the proof into two subcases: {\bf case $b_1$: $\frac{y}{\alpha}\geq\frac{1}{2}$} and
{\bf case $b_2: \frac{y}{\alpha}\in(0,\frac{1}{2})$.
}

{\bf case $b_1$: $\frac{y}{\alpha}\geq\frac{1}{2}$}. By Lemmas \ref{Lemma37}, \ref{Lemma23} and \ref{Lemma24}, we have
\begin{equation}\aligned\label{YYY1}
-\frac{\partial}{\partial x}\mathcal{W}_{\frac{1}{2\pi}}(\alpha;z)&=\mathcal{C}(\alpha,x,y)
\cdot
\Big(
\pi \alpha^2\cdot\big(1+\sum_{n=2}^\infty n^3e^{-\alpha\pi y (n^2-1)}\cdot\frac{\vartheta_Y(\frac{y}{\alpha};nx)}{\vartheta_Y(\frac{y}{\alpha};x)}\\
&\;\;+\frac{\vartheta_{XY}(\frac{y}{\alpha};x)}{\vartheta_Y(\frac{y}{\alpha};x)}
+\sum_{n=2}^\infty ne^{-\alpha\pi y (n^2-1)}\frac{\vartheta_{XY}(\frac{y}{\alpha};nx)}{\vartheta_{Y}(\frac{y}{\alpha};x)}
\Big)\\
&\geq
\pi\cdot\mathcal{C}(\alpha,x,y)
\cdot
\Big(
 \alpha^2\cdot\big(1-\sum_{n=2}^\infty n^4e^{-\alpha\pi y (n^2-1)}\cdot\frac{1+\mu(\frac{y}{\alpha})}{1-\mu(\frac{y}{\alpha})}\big)\\
&\;\;- \frac{1+\nu(\frac{y}{\alpha})}{1+\mu(\frac{y}{\alpha})}
-\sum_{n=2}^\infty n^2e^{-\alpha\pi y (n^2-1)}\cdot\frac{1+\nu(\frac{y}{\alpha})}{1-\mu(\frac{y}{\alpha})}
\Big)\\
\endaligned\end{equation}
Since $\frac{y}{\alpha}\geq\frac{1}{2}$, and $\mu, \nu$
are decreasing functions by \eqref{mmmx}, then
\begin{equation}\aligned\label{P1}
\frac{1+\nu(\frac{y}{\alpha})}{1-\mu(\frac{y}{\alpha})}
\leq\frac{1+\nu(\frac{1}{2})}{1-\mu(\frac{1}{2})}=1.186694067\cdots,\;
\frac{1+\mu(\frac{y}{\alpha})}{1-\mu(\frac{y}{\alpha})}
\leq\frac{1+\nu(\frac{1}{2})}{1-\mu(\frac{1}{2})}=1.074612508\cdots.
\endaligned\end{equation}
For $\frac{1+\nu(x)}{1+\mu(x)}$, there still holds that
\begin{equation}\aligned\label{H100}
\frac{1+\nu(x)}{1+\mu(x)}\;\;\hbox{is decreasing with}\;\;x\geq\frac{1}{2}.
\endaligned\end{equation}
The fact \eqref{H100} can be checked directly by taking derivative in view of \eqref{mmmx}, the details is omitted here.
Therefore, for $\frac{y}{\alpha}\geq\frac{1}{2}$, it holds that
\begin{equation}\aligned\label{P2}
\frac{1+\nu(x)}{1+\mu(x)}\leq\frac{1+\nu(\frac{1}{2})}{1+\mu(\frac{1}{2})}=1.104299511\cdots.
\endaligned\end{equation}
Denote the error terms in \eqref{YYY1} by
\begin{equation}\aligned\nonumber
\sigma_1:=\sum_{n=2}^\infty n^4e^{-\alpha\pi y (n^2-1)}\cdot\frac{1+\mu(\frac{y}{\alpha})}{1-\mu(\frac{y}{\alpha})}),\;
\sigma_2:=\sum_{n=2}^\infty n^2e^{-\alpha\pi y (n^2-1)}\cdot\frac{1+\nu(\frac{y}{\alpha})}{1-\mu(\frac{y}{\alpha})}
.
\endaligned\end{equation}
Then by \eqref{P1},
\begin{equation}\aligned\label{P3}
\sigma_1&\leq \frac{1+\mu(\frac{1}{2})}{1-\mu(\frac{1}{2})}\sum_{n=2}^\infty n^4e^{-\alpha\pi y (n^2-1)}
\leq\frac{1+\mu(\frac{1}{2})}{1-\mu(\frac{1}{2})}\sum_{n=2}^\infty n^4e^{-1.1\cdot\pi\cdot \frac{\sqrt3}{2} (n^2-1)}
\leq2.169\cdot 10^{-3},\\
\sigma_2&\leq\frac{1+\nu(\frac{1}{2})}{1-\mu(\frac{1}{2})}\sum_{n=2}^\infty n^2e^{-\alpha\pi y (n^2-1)}\leq \frac{1+\nu(\frac{1}{2})}{1-\mu(\frac{1}{2})}\sum_{n=2}^\infty n^2e^{-1.1\cdot\pi\cdot \frac{\sqrt3}{2} (n^2-1)}
\leq6.75\cdot 10^{-4}.
\endaligned\end{equation}

Therefore by \eqref{YYY1}, \eqref{P2} and \eqref{P3},
\begin{equation}\aligned\label{YYY1}
-\frac{\partial}{\partial x}\mathcal{W}_{\frac{1}{2\pi}}(\alpha;z)
&\geq
\pi\cdot\mathcal{C}(\alpha,x,y)
\cdot
\Big(
 \alpha^2\cdot\big(1-\sigma_1\big)- \frac{1+\nu(\frac{y}{\alpha})}{1+\mu(\frac{y}{\alpha})}
-\sigma_2
\Big)\\
&\geq
\pi\cdot\mathcal{C}(\alpha,x,y)
\cdot
\Big(1.1^2\cdot(1-2.169\cdot 10^{-3})-1.105-6.75\cdot 10^{-4}\Big)\\
&\geq
\pi\cdot\mathcal{C}(\alpha,x,y)
\cdot0.1017005100>0.
\endaligned\end{equation}
This completes the proof of {\bf case $b_1$}.

 It remains to prove {\bf case $b_2$}.

{\bf case $b_2: \frac{y}{\alpha}\in(0,\frac{1}{2})$.
}
By Lemmas \ref{Lemma37}, \ref{Lemma23} and \ref{Lemma25}, we have
\begin{equation}\aligned\label{YYYb}
-\frac{\partial}{\partial x}\mathcal{W}_{\frac{1}{2\pi}}(\alpha;z)&=\mathcal{C}(\alpha,x,y)
\cdot
\Big(
\pi \alpha^2\cdot\big(1+\sum_{n=2}^\infty n^3e^{-\alpha\pi y (n^2-1)}\cdot\frac{\vartheta_Y(\frac{y}{\alpha};nx)}{\vartheta_Y(\frac{y}{\alpha};x)}\\
&\;\;+\frac{\vartheta_{XY}(\frac{y}{\alpha};x)}{\vartheta_Y(\frac{y}{\alpha};x)}
+\sum_{n=2}^\infty ne^{-\alpha\pi y (n^2-1)}\frac{\vartheta_{XY}(\frac{y}{\alpha};nx)}{\vartheta_{Y}(\frac{y}{\alpha};x)}
\Big)\\
&\geq
\mathcal{C}(\alpha,x,y)
\cdot
\Big(
\pi \alpha^2\cdot(1-\frac{1}{\pi}\sum_{n=2}^\infty n^4 e^{-\alpha\pi \big((n^2-1)y-\frac{1}{4y}\big)})\\
&\;\;-\frac{3}{2}\frac{\alpha}{y}(1+\frac{\pi}{6}\frac{\alpha}{y})
-\frac{3}{2\pi}\frac{\alpha}{y}(1+\frac{\pi}{6}\frac{\alpha}{y})\cdot\sum_{n=2}^\infty n^2 e^{-\alpha\pi \big((n^2-1)y-\frac{1}{4y}\big)}
\Big).
\endaligned\end{equation}
Since $\frac{y}{\alpha}\in(0,\frac{1}{2})$, then $\frac{\alpha}{y}>2$ and $\alpha>2y\geq\sqrt3.$

Denote that the error terms in \eqref{YYYb}
\begin{equation}\aligned\label{P4}
\sigma_3:=\frac{1}{\pi}\sum_{n=2}^\infty n^4 e^{-\alpha\pi \big((n^2-1)y-\frac{1}{4y}\big)},\;\;
\sigma_4:=\frac{3}{2\pi}\frac{\alpha}{y}(1+\frac{\pi}{6}\frac{\alpha}{y})\cdot\sum_{n=2}^\infty n^2 e^{-\alpha\pi \big((n^2-1)y-\frac{1}{4y}\big)})
.
\endaligned\end{equation}
Then
\begin{equation}\aligned\label{P5}
\sigma_3&\leq\frac{1}{\pi}\sum_{n=2}^\infty n^4 e^{-\sqrt3\pi \big((n^2-1)\frac{\sqrt3}{2}-\frac{1}{2\sqrt3}\big)}\leq1.777\cdot10^{-6},\\
\sigma_4&\leq\frac{3}{\pi}(1+\frac{\pi}{3})\cdot\sum_{n=2}^\infty n^2 e^{-\sqrt3\pi \big((n^2-1)\frac{\sqrt3}{2}-\frac{1}{2\sqrt3}\big)})\leq2.727\cdot 10^{-5}
.
\endaligned\end{equation}
By \eqref{YYYb}, \eqref{P4} and \eqref{P5}, we have
\begin{equation}\aligned\label{YYYb}
-\frac{\partial}{\partial x}\mathcal{W}_{\frac{1}{2\pi}}(\alpha;z)
&\geq
\mathcal{C}(\alpha,x,y)
\cdot
\Big(
\pi \alpha^2\cdot(1-\sigma_3)-(\pi+3)
-\sigma_4
\Big)\\
&\geq
\mathcal{C}(\alpha,x,y)
\cdot
\Big(
3\pi \cdot(1-\sigma_3)-(\pi+3)
-\sigma_4
\Big)\\
&>0.
\endaligned\end{equation}
This proves {\bf case $b_2: \frac{y}{\alpha}\in(0,\frac{1}{2})$.
} The proof is complete.

\end{proof}

Next we are going to prove
{\bf case a} of Theorem \ref{Th32}. Namely,
\begin{proposition}\label{Th32b} Assume that $\alpha\in(1,1.1]$. Then for $b=\frac{1}{2\pi}$,
\begin{equation}\aligned
\frac{\partial}{\partial x}\mathcal{W}_{\frac{1}{2\pi}}(\alpha;z)<0,\;\;\hbox{for}\;\;z\in\mathcal{D}_{\mathcal{G}}.
\endaligned\end{equation}
\end{proposition}

We start with a lemma which is direct consequence of Lemma \ref{Lemma36}.
\begin{lemma} \label{Lemma38}
\begin{equation}\aligned\nonumber
-\frac{\partial}{\partial x}\mathcal{W}_{\frac{1}{2\pi}}(\alpha;z)=&8\pi\alpha^{-\frac{5}{2}}y^{\frac{3}{2}}\cdot\Big(
\sum_{n=1}^\infty\sum_{m=1}^\infty
n^3m\big(
\alpha^2 e^{-\pi y(\alpha n^2+\frac{1}{\alpha}m^2)}-e^{-\pi y(\alpha m^2+\frac{1}{\alpha}n^2)}
\big)\cdot \sin(2mn\pi x)
\Big).
\endaligned\end{equation}
\end{lemma}

Given by Lemma \ref{Lemma38}, we denote for convenience that
\begin{equation}\aligned\label{A1}
\mathcal{A}_{n,m}(\alpha;y):=n^3m\cdot\big(
\alpha^2 e^{-\pi y(\alpha n^2+\frac{1}{\alpha}m^2)}-e^{-\pi y(\alpha m^2+\frac{1}{\alpha}n^2)})
.
\endaligned\end{equation}
Then by Lemma \ref{Lemma38}, one rewrites $-\frac{\partial}{\partial x}\mathcal{W}_{\frac{1}{2\pi}}(\alpha;z)$
by
\begin{equation}\aligned\nonumber
-\frac{\partial}{\partial x}\mathcal{W}_{\frac{1}{2\pi}}(\alpha;z)=&8\pi\alpha^{-\frac{5}{2}}y^{\frac{3}{2}}\cdot\Big(
\sum_{n=1}^\infty\sum_{m=1}^\infty
\mathcal{A}_{n,m}(\alpha;y)\cdot \sin(2mn\pi x)
\Big).
\endaligned\end{equation}
Therefore, Proposition \ref{Th32b} is equivalent to
\begin{lemma}\label{Lemma39}
Assume that $\alpha\in(1,1.1]$. Then
\begin{equation}\aligned\label{D1}
\sum_{n=1}^\infty\sum_{m=1}^\infty
\mathcal{A}_{n,m}(\alpha;y)\cdot \sin(2mn\pi x)
>0
\endaligned\end{equation}
for $z=x+iy\in\mathcal{D}_{\mathcal{G}}$. Here $\mathcal{A}_{n,m}(\alpha;y)$ is defined in \eqref{A1}. In fact, we show that
\begin{equation}\aligned\nonumber
\sum_{n=1}^\infty\sum_{m=1}^\infty
\mathcal{A}_{n,m}(\alpha;y)\cdot \sin(2mn\pi x)
\geq\frac{1}{2}(\alpha^2-1)\sin(2\pi x)>0\;\;\hbox{for}\;\;z=x+iy\in\mathcal{D}_{\mathcal{G}}.
\endaligned\end{equation}

\end{lemma}

To prove Lemma \ref{Lemma39}, we split the double infinite sum \eqref{D1} into three parts as follows:
\begin{equation}\aligned\label{D2}
\sum_{n=1}^\infty\sum_{m=1}^\infty
=\sum_{m=n}+\sum_{m=1}^\infty\sum_{n=m+1}^\infty+\sum_{n=1}^\infty\sum_{m=n+1}^\infty.
\endaligned\end{equation}
To estimate each part in \eqref{D2}, we establish the following two lemmas.

\begin{lemma}\label{Lemma310}
\begin{equation}\aligned\label{D3}
\mid\sum_{n=m+1}^\infty
\mathcal{A}_{n,m}(\alpha;y)\cdot \sin(2mn\pi x)\mid
\leq B\cdot \mathcal{A}_{m,m}(\alpha;y)\cdot \mid\sin(2m^2\pi x)\mid,
\endaligned\end{equation}
\end{lemma}
where the constant $B$ is defined by
\begin{equation}\aligned\label{B100}
B:=\frac{2^6\alpha_0\pi y e^{-3\pi y\alpha_0}}
{1-2^6e^{-5\pi y\alpha_0}},
\endaligned\end{equation}
where $\alpha_0$ is an constant belonging to $(\frac{1}{\alpha},\alpha)$.
\begin{proof}
\begin{equation}\aligned\label{D3a}
\mid\sum_{n=m+1}^\infty
\frac{\mathcal{A}_{n,m}(\alpha;y)\cdot \sin(2mn\pi x)}{\mathcal{A}_{m,m}(\alpha;y)\cdot \sin(2m^2\pi x)}
\mid
&\leq\sum_{n=m+1}^\infty(\frac{n}{m})^4
\mid
\alpha\frac{\alpha e^{-\pi y\alpha(n^2-m^2)}-\frac{1}{\alpha}e^{-\pi y\frac{1}{\alpha}(n^2-m^2)}}{\alpha^2-1}
\mid\\
&=\sum_{n=m+1}^\infty(\frac{n}{m})^4\big(\alpha_0\pi y(n^2-m^2)-1 \big) e^{-\pi y \alpha_0(n^2-m^2)},\;\;\alpha_0\in(\frac{1}{\alpha},\alpha)\\
&\leq\sum_{n=m+1}^\infty(\frac{n}{m})^4\alpha_0\pi y n^2 e^{-\pi y \alpha_0(n^2-m^2)}.
\endaligned\end{equation}
Here we used the mean value Theorem to estimate
\begin{equation}\label{BB1}
\frac{\alpha e^{-\pi y\alpha(n^2-m^2)}-\frac{1}{\alpha}e^{-\pi y\frac{1}{\alpha}(n^2-m^2)}}{\alpha^2-1}.
\end{equation}
Continuing by \eqref{D3a}, we deform by letting $k=n-m$ that
\begin{equation}\aligned\label{D3b}
\sum_{n=m+1}^\infty(\frac{n}{m})^4\alpha_0\pi y n^2 e^{-\pi y \alpha_0(n^2-m^2)}
=\sum_{k=1}^\infty(\frac{m+k}{m})^4(m+k)^2\alpha_0\pi y e^{-\alpha_0\pi y(2m+k)k}.
\endaligned\end{equation}
To simplify the notations, we denote that
\begin{equation}\aligned\label{D3c}
b_k:=(\frac{m+k}{m})^4(m+k)^2\alpha_0\pi y e^{-\alpha_0\pi y(2m+k)k}, m\geq1.
\endaligned\end{equation}
Then, by \eqref{D3a}, \eqref{D3b} and \eqref{D3c}, one has
\begin{equation}\aligned\label{D3d}
\mid\sum_{n=m+1}^\infty
\frac{\mathcal{A}_{n,m}(\alpha;y)\cdot \sin(2mn\pi x)}{\mathcal{A}_{m,m}(\alpha;y)\cdot \sin(2m^2\pi x)}
\mid
\leq\sum_{k=1}^\infty b_k.
\endaligned\end{equation}
To provide an upper bound of $\sum_{k=1}^\infty b_k$, we estimate that
\begin{equation}\aligned\label{D3e}
\frac{b_{k+1}}{b_k}&=(1+\frac{1}{m+k})^6 e^{-\pi y\alpha_0(1+2k+2m)}\\
&\leq 2^6 e^{-5\pi y\alpha_0}:=q.
\endaligned\end{equation}
By \eqref{D3e}, $\sum_{k=1}^\infty b_k$ is controlled by a geometric sequence and then
\begin{equation}\aligned\label{D3f}
\sum_{k=1}^\infty b_k
\leq\frac{b_1}{1-q}\leq\frac{2^6\alpha_0\pi y e^{-3\pi y\alpha_0}}
{1-2^6e^{-5\pi y\alpha_0}}.
\endaligned\end{equation}
The desired result then follows by \eqref{D3d} and \eqref{D3f}.

\end{proof}

Dual to Lemma \ref{Lemma310}, we have
\begin{lemma}\label{Lemma311}
\begin{equation}\aligned\nonumber
\mid\sum_{m=n+1}^\infty
\mathcal{A}_{n,m}(\alpha;y)\cdot \sin(2mn\pi x)\mid
\leq B\cdot \mathcal{A}_{n,n}(\alpha;y)\cdot \mid\sin(2n^2\pi x)\mid,
\endaligned\end{equation}
where the constant $B$ is defined in \eqref{B100}.
\end{lemma}
\begin{proof} The proof is similar to that of Lemma \ref{Lemma310}.
\begin{equation}\aligned\label{F3a}
\mid\sum_{m=n+1}^\infty
\frac{\mathcal{A}_{n,m}(\alpha;y)\cdot \sin(2mn\pi x)}{\mathcal{A}_{m,m}(\alpha;y)\cdot \sin(2m^2\pi x)}
\mid
&\leq\sum_{m=n+1}^\infty(\frac{m}{n})^4
\mid
\alpha\frac{\alpha e^{-\pi y\frac{1}{\alpha}(m^2-n^2)}-\frac{1}{\alpha}e^{-\pi y{\alpha}(m^2-n^2)}}{\alpha^2-1}
\mid\\
&=\sum_{n=m+1}^\infty(\frac{n}{m})^4\big(\alpha_0\pi y(m^2-n^2)+1 \big) e^{-\pi y \alpha_0(n^2-m^2)},\;\;\alpha_0\in(\frac{1}{\alpha},\alpha)\\
&\leq\sum_{n=m+1}^\infty(\frac{n}{m})^4\alpha_0\pi y m^2 e^{-\pi y \alpha_0(n^2-m^2)}.
\endaligned\end{equation}
Here we used the mean value Theorem for
\begin{equation}\aligned\label{F3b}
\frac{\alpha e^{-\pi y\frac{1}{\alpha}(m^2-n^2)}-\frac{1}{\alpha}e^{-\pi y{\alpha}(m^2-n^2)}}{\alpha^2-1},
\endaligned\end{equation}
which is a litter different from \eqref{BB1}.

Given by \eqref{F3a}, the rest of the proof is the same to the proof of Lemma \ref{Lemma310} by exchanging $m$ and $n$.

\end{proof}

We are ready to prove Lemma \ref{Lemma39}.
\begin{proof}{\bf Proof of Lemma \ref{Lemma39}.} As the splitting in \eqref{D2},
\begin{equation}\aligned\label{Ga1}
&\sum_{n=1}^\infty\sum_{m=1}^\infty
\mathcal{A}_{n,m}(\alpha;y)\cdot \sin(2mn\pi x)
=\sum_{k=1}^\infty
\mathcal{A}_{k,k}(\alpha;y)\cdot \sin(2k^2\pi x)\\
&+\sum_{m=1}^\infty\sum_{n=m+1}^\infty\mathcal{A}_{n,m}(\alpha;y)\cdot \sin(2mn\pi x)
+\sum_{n=1}^\infty\sum_{m=n+1}^\infty\mathcal{A}_{n,m}(\alpha;y)\cdot \sin(2mn\pi x)\\
&\geq\sum_{k=1}^\infty
\mathcal{A}_{k,k}(\alpha;y)\cdot \sin(2k^2\pi x)
-2B\sum_{m=1}^\infty\mathcal{A}_{m,m}(\alpha;y)\cdot \mid\sin(2m^2\pi x)\mid,
\endaligned\end{equation}
where the constant $B$ is defined in \eqref{B100}.
Note that $z\in\mathcal{D}_{\mathcal{G}}$ implies that $x\in(0,\frac{1}{2})$ and $y>\frac{\sqrt3}{2}$.
Then $\sin(2\pi x)>0$.
Therefore, by \eqref{Ga1},
\begin{equation}\aligned\label{Gaa2}
\sum_{n=1}^\infty&\sum_{m=1}^\infty
\mathcal{A}_{n,m}(\alpha;y)\cdot \sin(2mn\pi x)
\geq(1-2B)\mathcal{A}_{1,1}(\alpha;y)\sin(2\pi x)-(1+2B)\sum_{n=2}^\infty\mathcal{A}_{n,n}(\alpha;y)\cdot \mid\sin(2n^2\pi x)\mid\\
&\geq
(1-2B)\mathcal{A}_{1,1}(\alpha;y)\sin(2\pi x)\cdot\Big(1-\frac{1+2B}{1-2B}\sum_{n=2}^\infty\mathcal{A}_{n,n}(\alpha;y)\cdot \mid\frac{\sin(2n^2\pi x)}{\sin(2\pi x)}\mid\Big)\\
&\geq
(1-2B)\mathcal{A}_{1,1}(\alpha;y)\sin(2\pi x)\cdot\Big(1-\frac{1+2B}{1-2B}\sum_{n=2}^\infty n^2\mathcal{A}_{n,n}(\alpha;y)\Big)
\endaligned\end{equation}
Note that by \eqref{A1},
\begin{equation}\aligned\label{Gaa3}
\mathcal{A}_{n,n}=(\alpha^2-1)n^4 e^{-n\pi y(\alpha+\frac{1}{\alpha})}.
\endaligned\end{equation}
Then
\begin{equation}\aligned\label{Gaa4}
\sum_{n=2}^\infty n^2\mathcal{A}_{n,n}=\sum_{n=2}^\infty n^6 e^{-n\pi y(\alpha+\frac{1}{\alpha})}\leq\sum_{n=2}^\infty n^6e^{-2n\pi y}
\leq\sum_{n=2}^\infty n^6e^{-n\sqrt3 \pi }\leq1.27\cdot 10^{-3}.
\endaligned\end{equation}
Therefore, by \eqref{Gaa2}, \eqref{Gaa3} and \eqref{Gaa4}, we have
\begin{equation}\aligned\label{Gaa5}
\sum_{n=1}^\infty\sum_{m=1}^\infty
\mathcal{A}_{n,m}(\alpha;y)\cdot \sin(2mn\pi x)
\geq(\alpha^2-1)\frac{9}{10}(1-2B)\sin(2\pi x)>0.
\endaligned\end{equation}
This completes the proof of Lemma \ref{Lemma39}.

\end{proof}

Therefore, Theorem \ref{Th32} is proved by Propositions \ref{Th32a} and \ref{Th32b}, and Proposition \ref{Th32b} is proved by Lemma \ref{Lemma39}.

\section{Analysis on the vertical line $\Gamma$}

In Theorem \ref{Th31}, we have established that
for $\alpha\geq1$,
\begin{equation}\aligned
\min_{z\in\mathbb{H}}\mathcal{W}_{\frac{1}{2\pi}}(\alpha;z)
=\min_{z\in\Gamma}\mathcal{W}_{\frac{1}{2\pi}}(\alpha;z),
\endaligned\end{equation}
where the vertical line $\Gamma$ is defined as
\begin{equation}\aligned\nonumber
\Gamma=\{
z\in\mathbb{H}: \Re(z)=\frac{1}{2},\; \Im(z)\geq\frac{\sqrt3}{2}
\}
\endaligned\end{equation}
see \eqref{Gfff}.

In this section, we aim to establish that

\begin{theorem}\label{Th41} Assume that $\alpha>1$. Then
\begin{equation}\aligned
\min_{z\in\Gamma}\mathcal{W}_{\frac{1}{2\pi}}(\alpha;z)\;\;\hbox{is achieved at}\;\; e^{i\frac{\pi}{3}}.
\endaligned\end{equation}
\end{theorem}

The proof of Theorem \ref{Th41} is based on the following Proposition

\begin{proposition}\label{Prop41}For $\alpha\geq$1 and $y\geq\frac{\sqrt3}{2}$,
\begin{equation}\aligned\nonumber
\partial_y\mathcal{W}_{\frac{1}{2\pi}}(\alpha;\frac{1}{2}+iy)\geq0.
\endaligned\end{equation}

\end{proposition}

We first state a lemma on zeros of $\partial_y\mathcal{W}_{\frac{1}{2\pi}}(\alpha;\frac{1}{2}+iy)$, which is deduced by Lemma \ref{Lemma35}
and Proposition 3.4(\cite{Bet2018}).

\begin{lemma}\label{Lemma41} Zeros of of $\partial_y\mathcal{W}_{\frac{1}{2\pi}}(\alpha;\frac{1}{2}+iy)$.
\begin{itemize}
  \item $\alpha=1$ is a first order zero of $\partial_y\mathcal{W}_{\frac{1}{2\pi}}(\alpha;\frac{1}{2}+iy)$ with respect to $\alpha$ for any $y>0$;
  \item $y=\frac{\sqrt3}{2}$ is a first order zero of $\partial_y\mathcal{W}_{\frac{1}{2\pi}}(\alpha;\frac{1}{2}+iy)$ with respect to $y$ for any $\alpha>0$.
\end{itemize}
Qualitatively,
\begin{equation}\aligned\label{HHH}
\lim_{a\rightarrow1, y\rightarrow\frac{\sqrt3}{2}}\frac{\partial_y\mathcal{W}_{\frac{1}{2\pi}}(\alpha;\frac{1}{2}+iy)}{(a-1)(y-\frac{\sqrt3}{2})}
=\partial_{yya}\mathcal{W}_{\frac{1}{2\pi}}(\alpha;\frac{1}{2}+iy)\mid_{a=1,y=\frac{\sqrt3}{2}}=1.127521373\cdots>0.
\endaligned\end{equation}
\end{lemma}

\begin{proof} The zero points property in first and second item are deduced by Lemma \ref{Lemma35} and Proposition 3.4 in \cite{Bet2018} respectively. \eqref{HHH} is computed by L'H\^opital's rule. The first order of the zeros is then followed by \eqref{HHH}.

\end{proof}

To prove Proposition \ref{Prop41}, based on Lemma \ref{Lemma41}, we divide its proof into four cases.
For convenience for stating the strategy, we denote that
\begin{equation}\aligned\nonumber
\mathcal{R}_a:&=\{(\alpha,y)\mid \alpha\in[1,1.2],\; y\in[\frac{\sqrt3}{2},1]
\},\\
\mathcal{R}_b:&=\{(\alpha,y)\mid \alpha\in[1,1.2],\; y\geq1
\},\\
\mathcal{R}_c:&=\{(\alpha,y)\mid \alpha\geq1.2,\; y\geq\frac{5}{6}\alpha
\},\\
\mathcal{R}_d:&=\{(\alpha,y)\mid \alpha\geq1.2,\; y\in[\frac{\sqrt3}{2},\frac{5}{6}\alpha]
\}.
\endaligned\end{equation}
Then
\begin{equation}\aligned\nonumber
\{(\alpha,y)\mid \alpha\geq1,\; y\geq\frac{\sqrt3}{2}
\}=\mathcal{R}_a\cup\mathcal{R}_b\cup\mathcal{R}_c\cup\mathcal{R}_d.
\endaligned\end{equation}

We shall prove that $\partial_y\mathcal{W}_{\frac{1}{2\pi}}(\alpha;\frac{1}{2}+iy)$ is nonnegative on $\mathcal{R}_a$,
$\mathcal{R}_b$, $\mathcal{R}_c$ and $\mathcal{R}_d$ respectively. In each region, we use different methods.
In Regions $\mathcal{R}_b$ and $\mathcal{R}_c$, we estimate directly of $\partial_y\mathcal{W}_{\frac{1}{2\pi}}(\alpha;\frac{1}{2}+iy)$
by its double sum and exponential expansion respectively. In Region $\mathcal{R}_d$, we estimate
$(\partial_{yy}+\frac{2}{y}\partial_y)\mathcal{W}_{\frac{1}{2\pi}}(\alpha;\frac{1}{2}+iy)$. While in Region $\mathcal{R}_a$
, we estimate $\partial_{yya}\mathcal{W}_{\frac{1}{2\pi}}(\alpha;\frac{1}{2}+iy)$. We prove the cases of
$\mathcal{R}_b$,
$\mathcal{R}_c$, $\mathcal{R}_d$ and $\mathcal{R}_a$ in the next four subsections respectively.

\subsection{Region $\mathcal{R}_b$: estimate of $\partial_y\mathcal{W}_{\frac{1}{2\pi}}(\alpha;\frac{1}{2}+iy)$}

\vskip0.3in

In this subsection, we shall prove that
\begin{lemma}\label{Lemma42} For $(\alpha,y)\in\mathcal{R}_b$, then $\partial_y\mathcal{W}_{\frac{1}{2\pi}}(\alpha;\frac{1}{2}+iy)\geq0$.
\end{lemma}

The proof of Lemma \ref{Lemma42} is based on the following Lemmas \ref{Lemma43} and \ref{Lemma44}.

\begin{lemma}\label{Lemma43} For $(\alpha,y)\in\mathcal{R}_b$, then

\begin{equation}\aligned\nonumber
\partial_y\mathcal{W}_{\frac{1}{2\pi}}(\alpha;\frac{1}{2}+iy)\geq
2\pi(\alpha^2-1) y^{\frac{1}{2}}e^{-\pi\frac{y}{\alpha}}
\mathcal{L}_b(\alpha;y),
\endaligned\end{equation}
where
\begin{equation}\aligned\nonumber
\mathcal{L}_b(\alpha;y):=
\frac{\pi y}{\alpha}-\frac{3}{2}-(\pi y\alpha-\frac{3}{2})\alpha^2 e^{-\pi y(\alpha-\frac{1}{\alpha})}
+(3\pi(1-B)-2(1+B)\frac{\alpha^2+1}{\alpha}y)(\alpha^2-1)e^{-\pi y\alpha}
.
\endaligned\end{equation}
The constant $B$ is very small is located in Lemma \ref{Lemma48}, i.e.,
\begin{equation}\aligned\nonumber
B=\frac{2^6\alpha_0\pi y e^{-3\pi y\alpha_0}}
{1-2^6e^{-5\pi y\alpha_0}},
\endaligned\end{equation}
where $\alpha_0\in(\frac{1}{\alpha},\alpha)$.
\end{lemma}

\begin{lemma}\label{Lemma44}
For $(\alpha,y)\in\mathcal{R}_b$, then
\begin{equation}\aligned\nonumber
\mathcal{L}_b(\alpha;y)\geq0.316(\alpha^2-1)\geq0.
\endaligned\end{equation}
\end{lemma}

\begin{proof} We first claim that $\partial_y\mathcal{L}_b(\alpha;y)>0$ for $(\alpha,y)\in\mathcal{R}_b$.
In fact,
\begin{equation}\aligned\nonumber
\partial_y\mathcal{L}_b(\alpha;y)
&=\pi(\frac{1}{\alpha}-\alpha^3 e^{-\pi y(\alpha-\frac{1}{\alpha})})
+\pi\alpha(\alpha^2-1)(\pi y\alpha-\frac{3}{2})e^{-\pi y(\alpha-\frac{1}{\alpha})}\\
&\geq\pi(\frac{1}{\alpha}-\alpha^3 e^{-\pi y(\alpha-\frac{1}{\alpha})}).
\endaligned\end{equation}
Now
\begin{equation}\aligned\nonumber
\frac{1}{\alpha}-\alpha^3 e^{-\pi y(\alpha-\frac{1}{\alpha})}
&=\alpha e^{-\frac{\pi y}{\alpha}}
\big(
\frac{1}{\alpha^2}e^{-\frac{\pi y}{\alpha}}-\alpha^2e^{-\pi y\alpha}
\big)\\
&= e^{-\frac{\pi y}{\alpha}}(\alpha^2-1)e^{-\pi y\alpha_b}(\pi y\alpha_b-2), \alpha_b\in(\frac{1}{\alpha},\alpha)
\endaligned\end{equation}
by mean value Theorem, which is positive since $\pi y\alpha_b-2\geq\frac{\pi}{\alpha}-2>0$.

Then it follows that
\begin{equation}\aligned\nonumber
\mathcal{L}_b(\alpha;y)\geq\mathcal{L}_b(\alpha;1)
=(\alpha^2-1)\cdot
\Big(
\frac{\frac{\pi }{\alpha}-\frac{3}{2}-(\pi \alpha-\frac{3}{2})\alpha^2 e^{-\pi (\alpha-\frac{1}{\alpha})}}{\alpha^2-1}\\
+(3\pi(1-B)-2(1+B)\frac{\alpha^2+1}{\alpha})e^{-\pi \alpha}
\Big)
\endaligned\end{equation}
The rest of the proof is based on elementary inequality
\begin{equation}\aligned\nonumber
\frac{\frac{\pi }{x}-\frac{3}{2}-(\pi x-\frac{3}{2})x^2 e^{-\pi (x-\frac{1}{x})}}{x^2-1}
+(3\pi(1-B)-2(1+B)\frac{x^2+1}{x})e^{-\pi x}
\geq0.316\cdots,\;\;\hbox{for}\;\;x\in[1,1.2].
\endaligned\end{equation}
Here $\frac{\frac{\pi }{x}-\frac{3}{2}-(\pi x-\frac{3}{2})x^2 e^{-\pi (x-\frac{1}{x})}}{x^2-1}$ has a removable singularity at $x=1$ and
\begin{equation}\aligned\nonumber
\lim_{x\rightarrow1}\frac{\frac{\pi }{x}-\frac{3}{2}-(\pi x-\frac{3}{2})x^2 e^{-\pi (x-\frac{1}{x})}}{x^2-1}
=\pi^2 - 3.5\pi + 1.5=0.374030114\cdots.
\endaligned\end{equation}
\end{proof}
It remains to prove Lemma \ref{Lemma43}. We use a deformation of $\partial_y\mathcal{W}_{\frac{1}{2\pi}}(\alpha;\frac{1}{2}+iy)$.

\begin{lemma}\label{Lemma45} A double sum expansion of $\partial_y\mathcal{W}_{\frac{1}{2\pi}}(\alpha;\frac{1}{2}+iy)$.
\begin{equation}\aligned\nonumber
\partial_y\mathcal{W}_{\frac{1}{2\pi}}(\alpha;\frac{1}{2}+iy)=
\frac{3}{2}y^{\frac{1}{2}}\Big(
-2\pi\sum_{n=1}^\infty n^2(e^{-\pi n^2\frac{y}{\alpha}}-\alpha^2 e^{-\pi n^2 y\alpha})\\
+4\pi\sum_{n=1}^\infty\sum_{m=1}^\infty(-1)^{mn}n^2
\big(
\alpha^2 e^{-\pi y(n^2\alpha+\frac{m^2}{\alpha})}-e^{-\pi y(m^2\alpha+\frac{n^2}{\alpha})}
\big)
\Big)\\
+y^{\frac{3}{2}}
\Big(
\frac{2\pi^2}{\alpha}\sum_{n=1}^\infty n^4(e^{-\pi n^2\frac{y}{\alpha}}-\alpha^4 e^{-\pi n^2 y\alpha})\\
+\frac{4\pi^2}{\alpha}\sum_{n=1}^\infty\sum_{m=1}^\infty(-1)^{mn}n^4
\big(
 e^{-\pi y(m^2\alpha+\frac{n^2}{\alpha})}-\alpha^4e^{-\pi y(n^2\alpha+\frac{m^2}{\alpha})}
\big)
\Big).
\endaligned\end{equation}

\end{lemma}

Lemma \ref{Lemma45} is based on the following Lemma

\begin{lemma}\label{Lemma46} The theta functions expression of $\partial_y\mathcal{W}_{\frac{1}{2\pi}}(\alpha;z)$
\begin{equation}\aligned\nonumber
\partial_y\mathcal{W}_{\frac{1}{2\pi}}(\alpha;z)
=\frac{3}{2}y^{\frac{1}{2}}
\Big(
\pi\alpha^2\sum_{n\in\mathbb{Z}}n^2 e^{-\alpha\pi y n^2}\vartheta(\frac{y}{\alpha};nx)
+\sum_{n\in\mathbb{Z}}e^{-\alpha\pi y n^2}\vartheta_X(\frac{y}{\alpha};nx)
\Big)\\
+y^{\frac{3}{2}}
\Big(
-\pi^2\alpha^3\sum_{n\in\mathbb{Z}}n^4 e^{-\alpha\pi y n^2}\vartheta(\frac{y}{\alpha};nx)
+\frac{1}{\alpha}\sum_{n\in\mathbb{Z}}e^{-\alpha\pi y n^2}\vartheta_{XX}(\frac{y}{\alpha};nx)
\Big).
\endaligned\end{equation}

\end{lemma}

Lemma \ref{Lemma46} is a direct consequence of Lemma \ref{Lemma32}.

Given by Lemma \ref{Lemma45}, We also need some auxiliary lemmas to prove Lemma \ref{Lemma43}.
\begin{lemma}\label{Lemma47} For $\alpha\in[1,7]$ and $y\geq\frac{\sqrt3}{2}$, it holds that
\begin{equation}\aligned\nonumber
2y^\frac{3}{2}\frac{\pi^2}{\alpha}
\sum_{n=2}^\infty n^4(e^{-\pi n^2\frac{y}{\alpha}}-\alpha^4 e^{-\pi n^2 y\alpha})
\geq3\pi y^\frac{1}{2}
\sum_{n=2}^\infty n^2(e^{-\pi n^2\frac{y}{\alpha}}-\alpha^2 e^{-\pi n^2 y\alpha})
\endaligned\end{equation}

\end{lemma}

\begin{proof} Denote that
\begin{equation}\aligned\label{Bn1}
\mathcal{B}_n(\alpha;y):=2y^\frac{3}{2}\frac{\pi^2}{\alpha}
n^4(e^{-\pi n^2\frac{y}{\alpha}}-\alpha^4 e^{-\pi n^2 y\alpha})
-3\pi y^\frac{1}{2}
 n^2(e^{-\pi n^2\frac{y}{\alpha}}-\alpha^2 e^{-\pi n^2 y\alpha}).
\endaligned\end{equation}
Then it equivalents to prove that
\begin{equation}\aligned\nonumber
\sum_{n=2}^\infty\mathcal{B}_n(\alpha;y)>0.
\endaligned\end{equation}
In fact, we shall show that
\begin{equation}\aligned\nonumber
\mathcal{B}_n(\alpha;y)>0,\;\;\hbox{for}\;\; n\geq2, \alpha\in[1,7]\;\;\hbox{and}\;\;y\geq\frac{\sqrt3}{2}.
\endaligned\end{equation}
By \eqref{Bn1}, one has
\begin{equation}\aligned\nonumber
\mathcal{B}_n(\alpha;y)
=n^2e^{-\pi n^2\frac{y}{\alpha}}\Big(
\frac{\pi}{\alpha} yn^2-\frac{3}{2}
-\alpha^2(\alpha\pi yn^2-\frac{3}{2})e^{-\pi yn^2(\alpha-\frac{1}{\alpha})}
\Big)
\endaligned\end{equation}
To show $\mathcal{B}_n(\alpha;y)$ is nonnegative in the desired region, it equivalents to show
$\frac{\pi}{\alpha} yn^2-\frac{3}{2}
-\alpha^2(\alpha\pi yn^2-\frac{3}{2})e^{-\pi yn^2(\alpha-\frac{1}{\alpha})}$ is nonnegative in the desired region. This is similar to the proof of Lemma \ref{Lemma44}. We consider a function
\begin{equation}\aligned\nonumber
\Big(
\frac{\pi}{\alpha} x-\frac{3}{2}
-\alpha^2(\alpha\pi x-\frac{3}{2})e^{-\pi x(\alpha-\frac{1}{\alpha})}
\Big), x=yn^2\geq2\sqrt3.
\endaligned\end{equation}

One has, since $yn^2\geq2\sqrt3$

\begin{equation}\aligned\nonumber
\Big(
\frac{\pi}{\alpha} -\frac{3}{2}
-\alpha^2(\alpha\pi yn^2-\frac{3}{2})e^{-\pi yn^2(\alpha-\frac{1}{\alpha})}
\Big)&\geq\Big(
2\sqrt3\frac{\pi}{\alpha} -\frac{3}{2}
-\alpha^2(\alpha2\sqrt3\pi -\frac{3}{2})e^{-2\sqrt3\pi (\alpha-\frac{1}{\alpha})}
\Big)\\
&\geq(\alpha^2-1)\frac{2\sqrt3\frac{\pi}{\alpha} -\frac{3}{2}
-\alpha^2(\alpha2\sqrt3\pi -\frac{3}{2})e^{-2\sqrt3\pi (\alpha-\frac{1}{\alpha})}
}{\alpha^2-1}\\
&\geq0.00113927433(\alpha^2-1)\;\;\hbox{for}\;\;\alpha\in[1,7].
\endaligned\end{equation}
Here $\frac{2\sqrt3\frac{\pi}{\alpha} -\frac{3}{2}
-\alpha^2(\alpha2\sqrt3\pi -\frac{3}{2})e^{-2\sqrt3\pi (\alpha-\frac{1}{\alpha})}
}{\alpha^2-1}$ has a removable singularity at $\alpha=1$ and
\begin{equation}\aligned\nonumber
\lim_{\alpha\rightarrow1}\frac{2\sqrt3\frac{\pi}{\alpha} -\frac{3}{2}
-\alpha^2(\alpha2\sqrt3\pi -\frac{3}{2})e^{-2\sqrt3\pi (\alpha-\frac{1}{\alpha})}
}{\alpha^2-1}=81.84546604\cdots.
\endaligned\end{equation}

\end{proof}

\begin{lemma}\label{Lemma48} For $\alpha\leq[1,1.2]$ and $y\geq\frac{\sqrt3}{2}$
\begin{equation}\aligned\nonumber
\sum_{n=1}^\infty\sum_{m=1}^\infty(-1)^{mn}n^4
\big(
 e^{-\pi y(m^2\alpha+\frac{n^2}{\alpha})}-\alpha^4e^{-\pi y(n^2\alpha+\frac{m^2}{\alpha})}
\big)
&\geq(1-B)(\alpha^4-1)e^{-\pi y(\alpha+\frac{1}{\alpha})},\\
\sum_{n=1}^\infty\sum_{m=1}^\infty(-1)^{mn}n^2
\big(
\alpha^2 e^{-\pi y(n^2\alpha+\frac{m^2}{\alpha})}-e^{-\pi y(m^2\alpha+\frac{n^2}{\alpha})}
\big)&\geq-(1+B)(\alpha^2-1)e^{-\pi y(\alpha+\frac{1}{\alpha})}.
\endaligned\end{equation}
Here the constant $B$ is defined in \eqref{B100} as
\begin{equation}\aligned\nonumber
B=\frac{2^6\alpha_0\pi y e^{-3\pi y\alpha_0}}
{1-2^6e^{-5\pi y\alpha_0}},
\endaligned\end{equation}
where $\alpha_0\in(\frac{1}{\alpha},\alpha)$.
\end{lemma}

The proof of Lemma \ref{Lemma48} is very similar to that of Lemmas \ref{Lemma310} and \ref{Lemma311}, hence we omit the details here.

\begin{proof}{\bf Proof of Lemma \ref{Lemma43}.} It follows by Lemmas \ref{Lemma45}, \ref{Lemma47} and \ref{Lemma48}.
\end{proof}
Therefore, the proof of Lemma \ref{Lemma42} is complete.

\subsection{Region $\mathcal{R}_c$: estimate of $\partial_y\mathcal{W}_{\frac{1}{2\pi}}(\alpha;\frac{1}{2}+iy)$}

\vskip0.3in

In this subsection, we shall prove that
\begin{lemma}\label{Lemma409} For $(\alpha,y)\in\mathcal{R}_c$, then $\partial_y\mathcal{W}_{\frac{1}{2\pi}}(\alpha;\frac{1}{2}+iy)>0$.
\end{lemma}

The proof of Lemma \ref{Lemma409} dues to the following Lemmas \ref{Lemma410}, \ref{Lemma411} and \ref{Lemma412}.
\begin{lemma}\label{Lemma410}For $(\alpha,y)\in\mathcal{R}_c$, then
\begin{equation}\aligned\nonumber
\partial_y\mathcal{W}_{\frac{1}{2\pi}}(\alpha;z)
\geq&\frac{3}{2}y^{\frac{1}{2}}\vartheta_X(\frac{y}{\alpha};0)
+y^{\frac{3}{2}}\frac{1}{\alpha}\vartheta_{XX}(\frac{y}{\alpha};0)
-2\pi y^{\frac{3}{2}}\alpha^3\vartheta(\frac{y}{\alpha};\frac{1}{2})(1+\epsilon_{c,3})e^{-\pi y}\\
&\;\;\;\;+2y^{\frac{3}{2}}\frac{1}{\alpha}(1+\epsilon_{c,4})\vartheta_{XX}(\frac{y}{\alpha};\frac{1}{2}) e^{-\alpha\pi y}.
\endaligned\end{equation}
Here $\epsilon_{c,3}$ and $\epsilon_{c,4}$ are very small and located in Lemmas \ref{Lemma413} and \ref{Lemma414} respectively.
\end{lemma}

The proof of \ref{Lemma410} is based on Lemmas \ref{Lemma46}, \ref{Lemma413} and \ref{Lemma414}.

We further simplify the lower bound of $\partial_y\mathcal{W}_{\frac{1}{2\pi}}(\alpha;z)$ in Lemma \ref{Lemma410}.
\begin{lemma}\label{Lemma411}For $(\alpha,y)\in\mathcal{R}_c$, then
\begin{equation}\aligned\nonumber
\partial_y\mathcal{W}_{\frac{1}{2\pi}}(\alpha;z)
\geq 2\pi y^{\frac{1}{2}} e^{-\frac{\pi y}{\alpha}}
\cdot
\Big(
\frac{\pi y}{\alpha}-\frac{3}{2}
-(1+\epsilon_{c,3}) y\alpha^3e^{-\pi y(\alpha-\frac{1}{\alpha})}
-2(1+\epsilon_{c,4})\frac{\pi y}{\alpha} e^{-\alpha\pi y}
\Big).
\endaligned\end{equation}
Here $\epsilon_{c,3}$ and $\epsilon_{c,4}$ are very small and located in Lemmas \ref{Lemma413} and \ref{Lemma414} respectively.
\end{lemma}

Now we can conclude that $\partial_y\mathcal{W}_{\frac{1}{2\pi}}(\alpha;z)$ is positive if $(\alpha,y)\in\mathcal{R}_c$
by the following
\begin{lemma}\label{Lemma412}For $(\alpha,y)\in\mathcal{R}_c$, then
\begin{equation}\aligned\nonumber
\frac{\pi y}{\alpha}-\frac{3}{2}
-(1+\epsilon_{c,3}) y\alpha^3e^{-\pi y(\alpha-\frac{1}{\alpha})}
-2(1+\epsilon_{c,4})\frac{\pi y}{\alpha} e^{-\alpha\pi y}
\geq\frac{1}{2}>0.
\endaligned\end{equation}
Here $\epsilon_{c,3}$ and $\epsilon_{c,4}$ are very small and located in Lemmas \ref{Lemma413} and \ref{Lemma414} respectively.
\end{lemma}

\begin{proof} Since $\frac{y}{\alpha}\geq\frac{5}{6}$, $\frac{\pi y}{\alpha}, -ye^{-\pi y(\alpha-\frac{1}{\alpha})}$,
$-ye^{-\alpha\pi y}$ are monotonically decreasing on $y$, one then has
\begin{equation}\aligned\nonumber
&\frac{\pi y}{\alpha}-\frac{3}{2}
-(1+\epsilon_{c,3}) y\alpha^3e^{-\pi y(\alpha-\frac{1}{\alpha})}
-2(1+\epsilon_{c,4})\frac{\pi y}{\alpha} e^{-\alpha\pi y}\\
\geq
&\frac{5\pi }{6}-\frac{3}{2}
-(1+\epsilon_{c,3}) \frac{5}{6}\alpha^4e^{-\frac{5\pi}{6} (\alpha^2-1)}
-2(1+\epsilon_{c,4})\frac{5\pi }{6} e^{-\frac{5\pi}{6}\alpha^2}.
\endaligned\end{equation}
The later is bigger than $\frac{1}{2}$ by a basic estimate and we omit the details here.

\end{proof}

The proof is based on the expression in Lemma \ref{Lemma46}. We then estimate each part in Lemma \ref{Lemma46} separately by Lemmas
\ref{Lemma413}-\ref{Lemma416}.

\begin{lemma}\label{Lemma413} Assume that $(\alpha,y)\in\mathcal{R}_c$, then
\begin{equation}\aligned\nonumber
\sum_{n\in\mathbb{Z}}n^2 e^{-\alpha\pi y n^2}\vartheta(\frac{y}{\alpha};\frac{n}{2})
&\geq 2 e^{-\pi\alpha y}\vartheta(\frac{y}{\alpha};\frac{1}{2})(1-\epsilon_{c,1}),\\
\sum_{n\in\mathbb{Z}}n^4 e^{-\alpha\pi y n^2}\vartheta(\frac{y}{\alpha};\frac{n}{2})
&\leq 2 e^{-\pi\alpha y}\vartheta(\frac{y}{\alpha};\frac{1}{2})(1+\epsilon_{c,3}).
\endaligned\end{equation}
Here
\begin{equation}\aligned\nonumber
\epsilon_{c,1}:
&=\frac{1+\sum_{k=1}^\infty e^{-\pi k^2\frac{y}{\alpha}}}{1-\sum_{k=1}^\infty e^{-\pi k^2\frac{y}{\alpha}}}\cdot
\sum_{n=2}^\infty n^2 e^{-\pi\alpha y(n^2-1)},\\
\epsilon_{c,3}:
&=\frac{1+\sum_{k=1}^\infty e^{-\pi k^2\frac{y}{\alpha}}}{1-\sum_{k=1}^\infty e^{-\pi k^2\frac{y}{\alpha}}}\cdot
\sum_{n=2}^\infty n^4 e^{-\pi\alpha y(n^2-1)}.
\endaligned\end{equation}
Numerically,
\begin{equation}\aligned\nonumber
\epsilon_{c,1}\leq5.68\cdot 10^{-4},\;\;
\epsilon_{c,3}\leq2.27\cdot 10^{-3}.
\endaligned\end{equation}

\end{lemma}

\begin{proof} We only prove the first one, the second one is much similar and we omit the details here.
We start with the deformation
\begin{equation}\aligned\label{K1a}
\sum_{n\in\mathbb{Z}}n^2 e^{-\alpha\pi y n^2}\vartheta(\frac{y}{\alpha};\frac{n}{2})
&=2\sum_{n=1}n^2 e^{-\alpha\pi y n^2}\vartheta(\frac{y}{\alpha};\frac{n}{2})\\
&= 2e^{-\alpha\pi y}\vartheta(\frac{y}{\alpha};\frac{1}{2})\cdot
\big(1+\sum_{n=2}^\infty n^2 e^{-\alpha\pi y(n^2-1)}\frac{\vartheta(\frac{y}{\alpha};\frac{n}{2})}{\vartheta(\frac{y}{\alpha};\frac{1}{2})}\big)
\endaligned\end{equation}
For $\frac{\vartheta(\frac{y}{\alpha};\frac{n}{2}}{\vartheta(\frac{y}{\alpha};\frac{1}{2})}$, one has
\begin{equation}\aligned\label{K1b}
\frac{\vartheta(\frac{y}{\alpha};\frac{n}{2}}{\vartheta(\frac{y}{\alpha};\frac{1}{2})}
\leq\frac{1+\sum_{k=1}^\infty e^{-\pi k^2\frac{y}{\alpha}}}{1-\sum_{k=1}^\infty e^{-\pi k^2\frac{y}{\alpha}}}
\endaligned\end{equation}
since for any $x$
\begin{equation}\aligned\nonumber
{1-\sum_{k=1}^\infty e^{-\pi k^2\frac{y}{\alpha}}}\leq{\vartheta(\frac{y}{\alpha};x})\leq
{1+\sum_{k=1}^\infty e^{-\pi k^2\frac{y}{\alpha}}}.
\endaligned\end{equation}
\eqref{K1a} and \eqref{K1b} yield the result.

\end{proof}

\begin{lemma}\label{Lemma414}If $\frac{y}{\alpha}\geq\frac{5}{6}$, then
\begin{equation}\aligned\nonumber
\sum_{n\in\mathbb{Z}}e^{-\alpha\pi y n^2}\vartheta_X(\frac{y}{\alpha};\frac{n}{2})
&\geq\vartheta_X(\frac{y}{\alpha};0)+ 2 e^{-\pi\alpha y}\vartheta_X(\frac{y}{\alpha};\frac{1}{2})(1-\epsilon_{c,2}),\\
\sum_{n\in\mathbb{Z}}e^{-\alpha\pi y n^2}\vartheta_{XX}(\frac{y}{\alpha};\frac{n}{2})
&\geq\vartheta_{XX}(\frac{y}{\alpha};0)+ 2 e^{-\pi\alpha y}\vartheta_{XX}(\frac{y}{\alpha};\frac{1}{2})(1+\epsilon_{c,4}).
\endaligned\end{equation}
Here
\begin{equation}\aligned\nonumber
\epsilon_{c,2}:
&=\frac{\sum_{k=1}^\infty k^2e^{-\pi (k^2-1)\frac{y}{\alpha}}}{1-4e^{-3\pi\frac{y}{\alpha}}}\cdot
\sum_{n=2}^\infty  e^{-\pi\alpha y(n^2-1)},\\
\epsilon_{c,4}:&=\frac{\sum_{k=1}^\infty k^4e^{-\pi (k^2-1)\frac{y}{\alpha}}}{1-16e^{-3\pi\frac{y}{\alpha}}}\cdot
\sum_{n=2}^\infty  e^{-\pi\alpha y(n^2-1)}.
\endaligned\end{equation}
Numerically,
\begin{equation}\aligned\nonumber
\epsilon_{c,2}\leq1.23\cdot 10^{-5},\;\;
\epsilon_{c,4}\leq1.24\cdot 10^{-5}.
\endaligned\end{equation}

\end{lemma}

\begin{proof} The proof the second one is very similar to the first one and then we only provide the proof of the first one here.

Deforming the expression, one has
\begin{equation}\aligned\label{K2a}
\sum_{n\in\mathbb{Z}}e^{-\alpha\pi y n^2}\vartheta_X(\frac{y}{\alpha};\frac{n}{2})
&=\vartheta_X(\frac{y}{\alpha};0)+2\sum_{n=1}^\infty e^{-\alpha\pi y n^2}\vartheta_X(\frac{y}{\alpha};\frac{n}{2})\\
&=\vartheta_X(\frac{y}{\alpha};0)+2 e^{-\frac{\pi y}{\alpha}}\vartheta_X(\frac{y}{\alpha};\frac{1}{2})\cdot
\big(
1+\sum_{n=2}^\infty e^{-\alpha\pi y(n^2-1)}\frac{\vartheta_X(\frac{y}{\alpha};\frac{n}{2})}{\vartheta_X(\frac{y}{\alpha};\frac{1}{2})}
\big).
\endaligned\end{equation}
For $\frac{\vartheta_X(\frac{y}{\alpha};\frac{n}{2})}{\vartheta_X(\frac{y}{\alpha};\frac{1}{2})}$, one has
\begin{equation}\aligned\label{K2b}
\mid\frac{\vartheta_X(\frac{y}{\alpha};\frac{n}{2})}{\vartheta_X(\frac{y}{\alpha};\frac{1}{2})}\mid
&=\frac
{\mid\sum_{k=1}^\infty(-1)^{kn}k^2 e^{-k^2\frac{\pi y}{\alpha}}\mid}
{\sum_{k=1}^\infty(-1)^{k-1}k^2 e^{-k^2\frac{\pi y}{\alpha}}}\\
&\leq\frac
{\sum_{k=1}^\infty k^2 e^{-(k^2-1)\frac{\pi y}{\alpha}}}
{1-\sum_{k=2}^\infty(-1)^{k-1}k^2 e^{-(k^2-1)\frac{\pi y}{\alpha}}}\\
&\leq\frac
{\sum_{k=1}^\infty k^2 e^{-(k^2-1)\frac{\pi y}{\alpha}}}
{1-4e^{-3\frac{\pi y}{\alpha}}}.
\endaligned\end{equation}
The result follows by \eqref{K2a} and \eqref{K2b}.

\end{proof}

\begin{lemma}\label{Lemma415} For $\frac{y}{\alpha}>\frac{3}{2\pi}$,
\begin{equation}\aligned\nonumber
\frac{3}{2}y^{\frac{1}{2}}\vartheta_X(\frac{y}{\alpha};0)
+y^{\frac{3}{2}}\frac{1}{\alpha}\vartheta_{XX}(\frac{y}{\alpha};0)
\geq2\pi y^{\frac{1}{2}}(\frac{\pi y}{\alpha}-\frac{3}{2}) e^{-\frac{\pi y}{\alpha}}.
\endaligned\end{equation}

\end{lemma}

\begin{proof} In using the explicit expression of $\vartheta_X, \vartheta_{XX}$, one has
\begin{equation}\aligned\nonumber
\frac{3}{2}y^{\frac{1}{2}}\vartheta_X(\frac{y}{\alpha};0)
+y^{\frac{3}{2}}\frac{1}{\alpha}\vartheta_{XX}(\frac{y}{\alpha};0)
=2\pi y^{\frac{1}{2}}\cdot
\Big(
(\frac{\pi y}{\alpha}-\frac{3}{2})e^{-\frac{\pi y}{\alpha}}
+\sum_{n=2}^\infty(\frac{\pi y}{\alpha}n^4-\frac{3}{2}n^2) e^{-\pi n^2\frac{y}{\alpha}}
\Big).
\endaligned\end{equation}
The result then follows.

\end{proof}

\begin{lemma}\label{Lemma416}
For $\frac{y}{\alpha}\geq\frac{5}{6}$,
\begin{equation}\aligned\nonumber
\vartheta(\frac{y}{\alpha};\frac{1}{2})\leq1,\;\;
\mid\vartheta_{XX}(\frac{y}{\alpha};\frac{1}{2})\mid\leq2\pi^2 e^{-\frac{\pi y}{\alpha}}.
\endaligned\end{equation}

\end{lemma}
\begin{proof} In using the explicit expression of $\vartheta, \vartheta_{XX}$, one has
\begin{equation}\aligned\nonumber
\vartheta(\frac{y}{\alpha};\frac{1}{2})&=1+2\sum_{n=1}^\infty(-1)^n e^{-n^2\frac{\pi y}{\alpha}}\\
\mid\vartheta_{XX}(\frac{y}{\alpha};\frac{1}{2})\mid&=2\pi^2\sum_{n=1}^\infty(-1)^{n-1} n^4e^{-n^2\frac{\pi y}{\alpha}}\\
&=2\pi^2 e^{-\frac{\pi y}{\alpha}}(1+\sum_{n=2}^\infty(-1)^{n-1} n^4e^{-(n^2-1)\frac{\pi y}{\alpha}})
\endaligned\end{equation}
The result then follows.

\end{proof}

\subsection{Region $\mathcal{R}_d$, estimate of
$(\partial_{yy}+\frac{2}{y}\partial_y)\mathcal{W}_{\frac{1}{2\pi}}(\alpha;\frac{1}{2}+iy)$}

In this subsection, we aim to prove that
\begin{lemma}\label{Lemmaw1} Assume that $(\alpha,y)\in\mathcal{R}_d$, then $(\partial_{yy}+\frac{2}{y}\partial_y)\mathcal{W}_{\frac{1}{2\pi}}(\alpha;\frac{1}{2}+iy)>0$.
\end{lemma}
We postpone the proof of Lemma \ref{Lemmaw1} and give the desired estimate we need as follows.

\begin{lemma}\label{Lemmaw2} Assume that $(\alpha,y)\in\mathcal{R}_d$, then $\partial_y\mathcal{W}_{\frac{1}{2\pi}}(\alpha;\frac{1}{2}+iy)\geq0$.
\end{lemma}

\begin{proof} Notice that
\begin{equation}\aligned\nonumber
\partial_{yy}+\frac{2}{y}\partial_y=y^{-2}\partial_y(y^2\partial_y).
\endaligned\end{equation}
Then by Lemma \ref{Lemmaw1}, one has
\begin{equation}\aligned\label{N2}
\partial_y(y^2\partial_y)\mathcal{W}_{\frac{1}{2\pi}}(\alpha;\frac{1}{2}+iy)>0\;\;\hbox{for}\;\;(\alpha,y)\in\mathcal{R}_d.
\endaligned\end{equation}

On the other hand, by Proposition 3.4 of B\'etermin \cite{Bet2018}, it holds that
\begin{equation}\aligned\label{aaF4}
(y^2\partial_y)\mathcal{W}_{\frac{1}{2\pi}}(\alpha;\frac{1}{2}+iy)\mid_{y=\frac{\sqrt3}{2}}=0\;\;\hbox{for}\;\;a>0.
\endaligned\end{equation}
\eqref{N2} and \eqref{aaF4} yield the result.

\end{proof}

In the rest of this subsection, we aim to prove Lemma \ref{Lemmaw1}.
We first have an identity for $\theta(\alpha;z)$, see \cite{Mon1988,LW2022}.
\begin{lemma}\label{Lemma419}It holds that
\begin{equation}\aligned\nonumber
(\partial_{yy}+\frac{2}{y}\partial_y)\theta(\alpha;z)
=&(\pi\alpha)^2\sum_{n,m}(n^2-\frac{(m+nx)^2}{y^2})^2 e^{-\pi\alpha(yn^2+\frac{(m+nx)^2}{y})}\\
&\;\;\;\;-\frac{2\pi\alpha}{y}\sum_{n,m}n^2 e^{-\pi\alpha(yn^2+\frac{(m+nx)^2}{y})}.
\endaligned\end{equation}

\end{lemma}

The following Lemma is deduced by Lemma \ref{Lemma419}.
\begin{lemma}\label{Lemma420} We have the differential identity for $\mathcal{W}_{\frac{1}{2\pi}}(\alpha;\frac{1}{2}+iy)$
\begin{equation}\aligned\nonumber
&(\partial_{yy}+\frac{2}{y}\partial_y)\mathcal{W}_{\frac{1}{2\pi}}(\alpha;\frac{1}{2}+iy)\\
=&(\pi\alpha)^2\sum_{n,m}(n^2-\frac{(m+\frac{n}{2})^2}{y^2})^2(yn^2+\frac{(m+\frac{n}{2})^2}{y}) e^{-\pi\alpha(yn^2+\frac{(m+\frac{n}{2})^2}{y})}\\
&\;\;\;\;+\frac{3}{y}\sum_{n,m}n^2 e^{-\pi\alpha(yn^2+\frac{(m+\frac{n}{2})^2}{y})}
-\frac{5}{2}\pi\alpha\sum_{n,m}(n^2-\frac{(m+\frac{n}{2})^2}{y^2})^2 e^{-\pi\alpha(yn^2+\frac{(m+\frac{n}{2})^2}{y})}\\
&\;\;\;\;-\frac{2\pi\alpha}{y}\sum_{n,m}n^2(yn^2+\frac{(m+\frac{n}{2})^2}{y}) e^{-\pi\alpha(yn^2+\frac{(m+\frac{n}{2})^2}{y})}
\endaligned\end{equation}
\end{lemma}

Proceed by Lemma \ref{Lemma420}, we deduce the lower bound of $(\partial_{yy}+\frac{2}{y}\partial_y)\mathcal{W}_{\frac{1}{2\pi}}(\alpha;\frac{1}{2}+iy)$.
\begin{lemma}[The lower bound of $(\partial_{yy}+\frac{2}{y}\partial_y)\mathcal{W}_{\frac{1}{2\pi}}(\alpha;\frac{1}{2}+iy)$]\label{Lemma421}
Assume that $(\alpha,y)\in\mathcal{R}_d$, then
\begin{equation}\aligned\nonumber
(\partial_{yy}+\frac{2}{y}\partial_y)\mathcal{W}_{\frac{1}{2\pi}}(\alpha;\frac{1}{2}+iy)
\geq
\pi\alpha y^{-4}e^{-\frac{\pi\alpha}{y}}\mathcal{L}_d(\alpha;y),
\endaligned\end{equation}
where
\begin{equation}\aligned\nonumber
\mathcal{L}_d(\alpha;y):=
\frac{2\pi\alpha}{y}-5(1+\epsilon_{d,1})+4\pi\alpha (y^2-\frac{1}{4})^2(y+\frac{1}{4y})e^{-\pi\alpha(y-\frac{3}{4y})}
-8(1+\epsilon_{d,2})y^3(y+\frac{1}{4y})e^{-\pi\alpha(y-\frac{3}{4y})}
\endaligned\end{equation}

\end{lemma}

Lemma \ref{Lemmaw1} is then proved by Lemma \ref{Lemma421} and following Lemma \ref{Lemma422}.
Lemma \ref{Lemma421} is proved by Lemma \ref{Lemma420} and Lemmas \ref{Lemma423}-\ref{Lemma426}.
\begin{lemma}[The positiveness of lower bound function in Lemma \ref{Lemma421}]\label{Lemma422}Assume that $(\alpha,y)\in\mathcal{R}_d$, then
\begin{equation}\aligned\nonumber
\mathcal{L}_d(\alpha;y)>0.
\endaligned\end{equation}

\end{lemma}

\begin{proof} We divide the proof into two cases, case a: $y\geq1$ and case b: $y\in[\frac{\sqrt3}{2},1]$.
For case a: $y\geq1$, then  $\big(4\pi\alpha(y^2-\frac{1}{4})^2-8(1+\epsilon_{d,2})y^3\big)\geq0$ since $\alpha\geq1.2$, then $\mathcal{L}_d(\alpha;y)>0$ immediately since $\frac{2\pi\alpha}{y}-5(1+\epsilon_{d,1})>0$.

For case b: $y\in[\frac{\sqrt3}{2},1]$, it is checked that $\frac{\partial}{\partial \alpha}\mathcal{L}_d(\alpha;y)>0$ for $\alpha\geq1.2$.
Then
\begin{equation}\aligned\label{dhhh1}
\mathcal{L}_d(\alpha;y)\geq&
\frac{2.4\pi }{y}-5(1+\epsilon_{d,1})+4.8 (y^2-\frac{1}{4})^2(y+\frac{1}{4y})e^{-1.2\pi(y-\frac{3}{4y})}\\
&\;\;\;\;-8(1+\epsilon_{d,2})y^3(y+\frac{1}{4y})e^{-1.2\pi(y-\frac{3}{4y})}.
\endaligned\end{equation}
The later explicit function in \eqref{dhhh1} has lower bound 7 for $y\in[\frac{\sqrt3}{2},1]$, the result then follows.

\end{proof}

In the following Lemmas \ref{Lemma423}-\ref{Lemma426}, we shall analyze each part of the identity in Lemma \ref{Lemma420}.
\begin{lemma}[A lower bound of double sum: first kind]\label{Lemma423}
\begin{equation}\aligned\nonumber
&\sum_{n,m}(n^2-\frac{(m+\frac{n}{2})^2}{y^2})^2(yn^2+\frac{(m+\frac{n}{2})^2}{y}) e^{-\pi\alpha(yn^2+\frac{(m+\frac{n}{2})^2}{y})}\\
\geq&\frac{2}{y^5}e^{-\frac{\pi\alpha}{y}}
+4(1-\frac{1}{4y^2})^2(y+\frac{1}{4y})e^{-\pi\alpha(y+\frac{1}{4y})}.
\endaligned\end{equation}
\end{lemma}

\begin{proof} The double sum evaluates at
$$(m,n)=\{(1,0),(-1,0)\}
\;\;\hbox{contributing}\;\;\frac{1}{y^5} e^{-\pi\frac{\alpha}{y}}\;\;\hbox{each}
$$
and
$$(m,n)=\{(0,1),(0,-1),(1,-1),(-1,1)\}
\;\;\hbox{contributing}\;\;(1-\frac{1}{4y^2})^2 (y+\frac{1}{4y})e^{-\pi\alpha(y+\frac{1}{4y})}\;\;\hbox{each}.
$$
The rest of other terms in the double sum all are positive and hence the result follows.
\end{proof}

\begin{lemma}[A lower bound of double sum: second kind]\label{Lemma424}
\begin{equation}\aligned\nonumber
\sum_{n,m}n^2 e^{-\pi\alpha(yn^2+\frac{(m+\frac{n}{2})^2}{y})}
\geq4e^{-\pi\alpha(y+\frac{1}{4y})}.
\endaligned\end{equation}
\end{lemma}

\begin{proof} The double sum can be evaluated at
$$(m,n)=\{(0,1),(0,-1),(1,-1),(-1,1)\}
\;\;\hbox{contributing}\;\; e^{-\pi\alpha(y+\frac{1}{4y})}\;\;\hbox{each}.
$$
The rest of other terms in the double sum all are positive and hence the result follows.

\end{proof}

\begin{lemma}[An upper bound of double sum: third kind]\label{Lemma425}
\begin{equation}\aligned\nonumber
\sum_{n,m}(n^2-\frac{(m+\frac{n}{2})^2}{y^2})^2 e^{-\pi\alpha(yn^2+\frac{(m+\frac{n}{2})^2}{y})}
\leq(1+\epsilon_{d,1})\frac{2}{y^4}e^{-\frac{\pi\alpha}{y}}+4(1-\frac{1}{4y^2})^2 e^{-\pi\alpha(y+\frac{1}{4y})},
\endaligned\end{equation}
where
\begin{equation}\aligned\nonumber
\epsilon_{d,1}
\leq 4y^4 e^{-\pi\alpha(4y-\frac{1}{y})}
+16y^4 e^{-4\pi\alpha y}\leq3.92 \cdot 10^{-4}.
\endaligned\end{equation}
\end{lemma}

\begin{proof} We deform the double sum as
\begin{equation}\aligned\label{dH1}
\sum_{n,m}(n^2-\frac{(m+\frac{n}{2})^2}{y^2})^2 e^{-\pi\alpha(yn^2+\frac{(m+\frac{n}{2})^2}{y})}
=\sum_{p\equiv q(\mod 2)}(p^2-\frac{q^2}{4y^2})^2 e^{-\pi\alpha(yp^2+\frac{q^2}{4y})}.
\endaligned\end{equation}
One then splits the double sum into four parts as
\begin{equation}\aligned\nonumber
(p,q)\in (a):
 p=\pm1, q=\pm1; (b): p=0, q=\pm2; (c): p=\pm2, q=0\;\;\hbox{and}\;\;(d): p\geq2, q\geq2.
\endaligned\end{equation}
Continuing with \eqref{dH1}, one has
\begin{equation}\aligned\label{dH2}
\sum_{n,m}(n^2-\frac{(m+\frac{n}{2})^2}{y^2})^2 e^{-\pi\alpha(yn^2+\frac{(m+\frac{n}{2})^2}{y})}
=&
4(1-\frac{1}{4y^2})^2 e^{-\pi\alpha(y+\frac{1}{4y})}
+\frac{2}{y^4}e^{-\frac{\pi\alpha}{y}}
+8 e^{-4\pi\alpha y}\\
+&\;\;\;\;\sum_{p\equiv q(\mod 2), p\geq2, q\geq2}(p^2-\frac{q^2}{4y^2})^2 e^{-\pi\alpha(yp^2+\frac{q^2}{4y})}.
\endaligned\end{equation}
The last term in \eqref{dH2} is very small and can be controlled by
\begin{equation}\aligned\label{dH3}
&\sum_{p\equiv q(\mod 2), p\geq2, q\geq2}(p^2-\frac{q^2}{4y^2})^2 e^{-\pi\alpha(yp^2+\frac{q^2}{4y})}\\
\leq&\sum_{p\equiv q(\mod 2), p\geq2, q\geq2}(p^4+\frac{q^4}{16y^4}) e^{-\pi\alpha(yp^2+\frac{q^2}{4y})}\\
\leq&\sum_{p\geq2}p^4 e^{-\pi\alpha y p^2}\sum_{q\geq2} e^{-\frac{\pi\alpha}{4y}q^2}
+\frac{1}{16y^4}\sum_{p\geq2} e^{-\pi\alpha y p^2}\sum_{q\geq2} q^4e^{-\frac{\pi\alpha}{4y}q^2}\\
\leq&
16 e^{-\pi\alpha(4y+\frac{1}{y})}\cdot d(\alpha;y).
\endaligned\end{equation}
Here $d(\alpha;y)$ is bounded by some constant and has the following expression
\begin{equation}\aligned\label{dH4}
d(\alpha;y):=
\sum_{p\geq2}(\frac{p}{2})^4 e^{-\pi\alpha y (p^2-4)}\sum_{q\geq2} e^{-\frac{\pi\alpha}{4y}(q^2-4)}
+\frac{1}{16y^4}\sum_{p\geq2} e^{-\pi\alpha y (p^2-4)}\sum_{q\geq2}(\frac{q}{2})^4 e^{-\frac{\pi\alpha}{4y}(q^2-4)}
.
\endaligned\end{equation}
Roughly, one has
\begin{equation}\aligned\label{dH4a}
d(\alpha;y)\leq2
.
\endaligned\end{equation}

Therefore, by \eqref{dH2}, \eqref{dH3} and \eqref{dH4a},
\begin{equation}\aligned\label{dH5}
\sum_{n,m}(n^2-\frac{(m+\frac{n}{2})^2}{y^2})^2 e^{-\pi\alpha(yn^2+\frac{(m+\frac{n}{2})^2}{y})}
\leq&
4(1-\frac{1}{4y^2})^2 e^{-\pi\alpha(y+\frac{1}{4y})}
+\frac{2}{y^4}e^{-\frac{\pi\alpha}{y}}\\
&\;\;\;\;+8 e^{-4\pi\alpha y}+32 e^{-\pi\alpha(4y+\frac{1}{y})}.
\endaligned\end{equation}
The proof is complete.

\end{proof}

\begin{lemma}[An upper bound of double sum: fourth kind]\label{Lemma426}
\begin{equation}\aligned\nonumber
\sum_{n,m}n^2(yn^2+\frac{(m+\frac{n}{2})^2}{y}) e^{-\pi\alpha(yn^2+\frac{(m+\frac{n}{2})^2}{y})}
\leq4(1+\epsilon_{d,2})(y+\frac{1}{4y})e^{-\pi\alpha(y+\frac{1}{4y})},
\endaligned\end{equation}
where
\begin{equation}\aligned\nonumber
\epsilon_{d,2}
\leq16e^{-3\pi\alpha y}(1+e^{-\frac{3\pi\alpha}{4y}})\leq9.27\cdot 10^{-4}.
\endaligned\end{equation}
\end{lemma}

\begin{proof}
We first deform the double sum as
\begin{equation}\aligned\label{dhh1}
\sum_{n,m}n^2(yn^2+\frac{(m+\frac{n}{2})^2}{y}) e^{-\pi\alpha(yn^2+\frac{(m+\frac{n}{2})^2}{y})}
=\sum_{p\equiv q(\mod2)}p^2(yp^2+\frac{q^2}{4y}) e^{-\pi\alpha(yp^2+\frac{q^2}{4y})}.
\endaligned\end{equation}
One then splits the double sum into four parts as
\begin{equation}\aligned\nonumber
(p,q)\in (a):
 p=\pm1, q=\pm1; (b): p=0, q=\pm2; (c): p=\pm2, q=0\;\;\hbox{and}\;\;(d): p\geq2, q\geq2.
\endaligned\end{equation}
Then by \eqref{dhh1},
\begin{equation}\aligned\label{dhh2}
\sum_{n,m}n^2(yn^2+\frac{(m+\frac{n}{2})^2}{y}) e^{-\pi\alpha(yn^2+\frac{(m+\frac{n}{2})^2}{y})}
=&
4(y+\frac{1}{4y})e^{-\pi\alpha(y+\frac{1}{4y})}+16y e^{-4\pi\alpha y}\\
&\;\;+\sum_{p\equiv q(\mod2), p\geq2, q\geq2}p^2(yp^2+\frac{q^2}{4y}) e^{-\pi\alpha(yp^2+\frac{q^2}{4y})}
\endaligned\end{equation}
The last term in \eqref{dhh2} is very small and can be controlled by
\end{proof}
\begin{equation}\aligned\label{dhh3}
&\sum_{p\equiv q(\mod2), p\geq2, q\geq2}p^2(yp^2+\frac{q^2}{4y}) e^{-\pi\alpha(yp^2+\frac{q^2}{4y})}\\
\leq&
y\sum_{p\geq2}p^4 e^{-\pi\alpha p^2}\sum_{q\geq2} e^{-\frac{\pi\alpha}{4y}q^2}
+\frac{1}{4y}\sum_{p\geq2}p^2 e^{-\pi\alpha p^2}\sum_{q\geq2} q^2e^{-\frac{\pi\alpha}{4y}q^2}.
\endaligned\end{equation}
The result then follows by \eqref{dhh2} and \eqref{dhh3} after some simple deformations and we omit the details here.

\subsection{Region $\mathcal{R}_a$, estimate of
$(\partial_{yy\alpha}+\frac{2}{y}\partial_{y\alpha})\mathcal{W}_{\frac{1}{2\pi}}(\alpha;\frac{1}{2}+iy)$}

In this section, we aim to establish that
\begin{lemma}\label{Lemmad1} Assume that $(\alpha,y)\in\mathcal{R}_a$, then
$(\partial_{yy\alpha}+\frac{2}{y}\partial_{y\alpha})\mathcal{W}_{\frac{1}{2\pi}}(\alpha;\frac{1}{2}+iy)>0$.

\end{lemma}

With Lemma \ref{Lemmad1}, one has
\begin{lemma}\label{Lemmad2} Assume that $(\alpha,y)\in\mathcal{R}_a$, then
$\partial_y\mathcal{W}_{\frac{1}{2\pi}}(\alpha;\frac{1}{2}+iy)\geq0$.

\end{lemma}

\begin{proof} Notice that
\begin{equation}\aligned\label{aF1}
\partial_{yy\alpha}+\frac{2}{y}\partial_{y\alpha}=\partial_\alpha (y^{-2}\partial_y(y^2\partial_y)).
\endaligned\end{equation}

By Lemma \ref{Lemma35}, one has
\begin{equation}\aligned\label{aF2}
y^{-2}\partial_y(y^2\partial_y)\mathcal{W}_{\frac{1}{2\pi}}(\alpha;\frac{1}{2}+iy)\mid_{a=1}=0\;\;\hbox{for}\;\;y>0.
\endaligned\end{equation}
Then by Lemma \ref{Lemmad1}, \eqref{aF1} and \eqref{aF2}
\begin{equation}\aligned\label{aF3}
\partial_y(y^2\partial_y)\mathcal{W}_{\frac{1}{2\pi}}(\alpha;\frac{1}{2}+iy)\geq0\;\;\hbox{for}\;\;(\alpha,y)\in\mathcal{R}_a.
\endaligned\end{equation}
On the other hand, by Proposition 3.4 of B\'etermin \cite{Bet2018}, it holds that
\begin{equation}\aligned\label{aF4}
(y^2\partial_y)\mathcal{W}_{\frac{1}{2\pi}}(\alpha;\frac{1}{2}+iy)\mid_{y=\frac{\sqrt3}{2}}=0\;\;\hbox{for}\;\;a>0.
\endaligned\end{equation}
Then by \eqref{aF3} and \eqref{aF4},
\begin{equation}\aligned\label{aF5}
(y^2\partial_y)\mathcal{W}_{\frac{1}{2\pi}}(\alpha;\frac{1}{2}+iy)\geq0\;\;\hbox{for}\;\;(\alpha,y)\in\mathcal{R}_a.
\endaligned\end{equation}
\eqref{aF5} yields the result.

\end{proof}

It remains to prove Lemma \ref{Lemmad1}.
We start from Lemma \ref{Lemma420}. After simple computation, one has
\begin{lemma}[An identity for $(\partial_{yy\alpha}+\frac{2}{y}\partial_{y\alpha})\mathcal{W}_{\frac{1}{2\pi}}(\alpha;\frac{1}{2}+iy)$]\label{Lemma429}
\begin{equation}\aligned\nonumber
(\partial_{yy\alpha}+\frac{2}{y}\partial_{y\alpha})\mathcal{W}_{\frac{1}{2\pi}}(\alpha;\frac{1}{2}+iy)
=&\frac{9}{2}\pi^2\alpha\sum_{n,m}(n^2-\frac{(m+\frac{n}{2})^2}{y^2})^2(yn^2+\frac{(m+\frac{n}{2})^2}{y}) e^{-\pi\alpha(yn^2+\frac{(m+\frac{n}{2})^2}{y})}\\
&\;\;+\frac{2\pi^2\alpha}{y}\sum_{n,m}n^2(yn^2+\frac{(m+\frac{n}{2})^2}{y})^2 e^{-\pi\alpha(yn^2+\frac{(m+\frac{n}{2})^2}{y})}\\
&\;\;-\frac{5\pi}{2}\sum_{n,m}(n^2-\frac{(m+\frac{n}{2})^2}{y^2})^2 e^{-\pi\alpha(yn^2+\frac{(m+\frac{n}{2})^2}{y})}\\
&\;\;-\frac{5\pi}{y}\sum_{n,m}n^2(yn^2+\frac{(m+\frac{n}{2})^2}{y}) e^{-\pi\alpha(yn^2+\frac{(m+\frac{n}{2})^2}{y})}\\
&\;\;-\pi^3\alpha^2\sum_{n,m}(n^2-\frac{(m+\frac{n}{2})^2}{y^2})^2(yn^2+\frac{(m+\frac{n}{2})^2}{y})^2 e^{-\pi\alpha(yn^2+\frac{(m+\frac{n}{2})^2}{y})}\\
\endaligned\end{equation}

\end{lemma}

Based on Lemma \ref{Lemma429}, we then state the following Lemma and postpone its proof to the late part of this subsection.

\begin{lemma}[A lower bound function of $(\partial_{yy\alpha}+\frac{2}{y}\partial_{y\alpha})\mathcal{W}_{\frac{1}{2\pi}}(\alpha;\frac{1}{2}+iy)$]\label{Lemma430}
Assume that $(\alpha,y)\in\mathcal{R}_a$, then
\begin{equation}\aligned\nonumber
(\partial_{yy\alpha}+\frac{2}{y}\partial_{y\alpha})\mathcal{W}_{\frac{1}{2\pi}}(\alpha;\frac{1}{2}+iy)
\geq\frac{\pi}{y^4}e^{-\frac{\pi\alpha}{y}}\cdot \mathcal{L}_a(\alpha;y).
\endaligned\end{equation}
Here
\begin{equation}\aligned\nonumber
\mathcal{L}_a(\alpha;y)=
\frac{9\pi\alpha}{y}-5-\frac{2\pi^2\alpha^2}{y^2}+
H(\alpha;y)e^{-\pi\alpha(y-\frac{3}{4y})},
\endaligned\end{equation}
and
\begin{equation}\aligned\nonumber
H(\alpha;y)=&
18\pi\alpha(y^2-\frac{1}{4})^2(y+\frac{1}{4y})+8\pi\alpha y^3(y+\frac{1}{4y})^2\\
&\;\;-10(y^2-\frac{1}{4})^2-20y^3(y+\frac{1}{4y})-4\pi^2\alpha^2(y^2-\frac{1}{4})^2(y+\frac{1}{4y})^2.
\endaligned\end{equation}

\end{lemma}

\begin{lemma}[The positiveness of the lower bound function in Lemma \ref{Lemma430}]\label{Lemma431}Assume that $(\alpha,y)\in\mathcal{R}_a$, then
\begin{equation}\aligned\nonumber
\mathcal{L}_a(\alpha;y)\geq\frac{1}{2}>0.
\endaligned\end{equation}
\end{lemma}

\begin{proof} Since $\mathcal{R}_a$ is a small finite region, we split it into 14 subregions to get the result.
\end{proof}

By Lemmas \ref{Lemma430} and \ref{Lemma431}, one gets Lemma \ref{Lemmad1}.

It remains to prove Lemma \ref{Lemma430}. We start from Lemma \ref{Lemma429}. There are five types double sum in Lemma \ref{Lemma429}, three of them are estimated in Lemmas \ref{Lemma423}, \ref{Lemma425}-\ref{Lemma426}.
We shall estimate the left two of them in the late part of this subsection. Lemma \ref{Lemma430} then follows from
Lemmas \ref{Lemma423}, \ref{Lemma425}-\ref{Lemma426} and \ref{Lemma432}-\ref{Lemma433}.

\begin{lemma}[A lower bound of double sum]\label{Lemma432}
\begin{equation}\aligned\nonumber
\sum_{n,m}n^2(yn^2+\frac{(m+\frac{n}{2})^2}{y})^2 e^{-\pi\alpha(yn^2+\frac{(m+\frac{n}{2})^2}{y})}
\geq4(y+\frac{1}{4y}))^2e^{-\pi\alpha(y+\frac{1}{4y})}.
\endaligned\end{equation}
\end{lemma}

\begin{proof} The double sum can be evaluated at
$$(m,n)=\{(0,1),(0,-1),(1,-1),(-1,1)\}
\;\;\hbox{contributing}\;\; (y+\frac{1}{4y}))^2e^{-\pi\alpha(y+\frac{1}{4y})}\;\;\hbox{each}.
$$
The rest of other terms in the double sum all are positive and hence the result follows.

\end{proof}

\begin{lemma}[An upper bound of double sum]\label{Lemma433}
\begin{equation}\aligned\nonumber
&\sum_{n,m}(n^2-\frac{(m+\frac{n}{2})^2}{y^2})^2(yn^2+\frac{(m+\frac{n}{2})^2}{y})^2 e^{-\pi\alpha(yn^2+\frac{(m+\frac{n}{2})^2}{y})}\\
\leq&\frac{2}{y^6}e^{-\frac{\pi\alpha}{y}}
+4(1-\frac{1}{4y^2})^2(y+\frac{1}{4y})^2e^{-\pi\alpha(y+\frac{1}{4y})}
+3\cdot 16^2 e^{-4\pi\alpha y}.
\endaligned\end{equation}

\end{lemma}

\begin{proof} We deform the double sum as
\begin{equation}\aligned\label{adH1}
&\sum_{n,m}(n^2-\frac{(m+\frac{n}{2})^2}{y^2})^2(yn^2+\frac{(m+\frac{n}{2})^2}{y})^2 e^{-\pi\alpha(yn^2+\frac{(m+\frac{n}{2})^2}{y})}\\
=&\sum_{p\equiv q(\mod 2)}(p^2-\frac{q^2}{4y^2})^2(yp^2+\frac{q^2}{4y})^2 e^{-\pi\alpha(yp^2+\frac{q^2}{4y})}.
\endaligned\end{equation}
One then splits the double sum into four parts as
\begin{equation}\aligned\nonumber
(p,q)\in (a):
 p=\pm1, q=\pm1; (b): p=0, q=\pm2; (c): p=\pm2, q=0\;\;\hbox{and}\;\;(d): p\geq2, q\geq2.
\endaligned\end{equation}
It follows that
\begin{equation}\aligned\label{adH2}
&\sum_{p\equiv q(\mod 2)}(p^2-\frac{q^2}{4y^2})^2(yp^2+\frac{q^2}{4y})^2 e^{-\pi\alpha(yp^2+\frac{q^2}{4y})}\\
=
&4(1-\frac{1}{4y^2})^2(y+\frac{1}{4y})^2 e^{-\pi\alpha(y+\frac{1}{4y})}
+\frac{2}{y^6}e^{-\frac{\pi\alpha}{y}}
+2\cdot 16^2 e^{-4\pi\alpha y}\\
+&\;\;\;\;\sum_{p\equiv q(\mod 2), p\geq2, q\geq2}(p^2-\frac{q^2}{4y^2})^2(yp^2+\frac{q^2}{4y})^2 e^{-\pi\alpha(yp^2+\frac{q^2}{4y})}.
\endaligned\end{equation}
The last term in \eqref{adH2} is very small and can be controlled by
\begin{equation}\aligned\label{adH3}
&\sum_{p\equiv q(\mod 2), p\geq2, q\geq2}(p^2-\frac{q^2}{4y^2})^2(yp^2+\frac{q^2}{4y})^2 e^{-\pi\alpha(yp^2+\frac{q^2}{4y})}\\
=&\sum_{p\equiv q(\mod 2), p\geq2, q\geq2}y^2(p^4-\frac{q^4}{16y^4})^2 e^{-\pi\alpha(yp^2+\frac{q^2}{4y})}\\
\leq&\sum_{p\equiv q(\mod 2), p\geq2, q\geq2}y^2(p^8+\frac{q^8}{16^2y^8}) e^{-\pi\alpha(yp^2+\frac{q^2}{4y})}\\
\leq&y^2\sum_{p\geq2}p^8 e^{-\pi\alpha y p^2}\sum_{q\geq2} e^{-\frac{\pi\alpha}{4y}q^2}
+\frac{1}{16^2y^6}\sum_{p\geq2} e^{-\pi\alpha y p^2}\sum_{q\geq2} q^8e^{-\frac{\pi\alpha}{4y}q^2}\\
\leq&
16^2 y^2 e^{-\pi\alpha(4y+\frac{1}{y})}\cdot d_2(\alpha;y)
\endaligned\end{equation}
Here $d_2(\alpha;y)$ is bounded by some constant and has the following expression
\begin{equation}\aligned\label{adH4}
d_2(\alpha;y):=
\sum_{p\geq2}(\frac{p}{2})^4 e^{-\pi\alpha y (p^2-4)}\sum_{q\geq2} e^{-\frac{\pi\alpha}{4y}(q^2-4)}
+\frac{1}{16y^4}\sum_{p\geq2} e^{-\pi\alpha y (p^2-4)}\sum_{q\geq2}(\frac{q}{2})^4 e^{-\frac{\pi\alpha}{4y}(q^2-4)}
.
\endaligned\end{equation}
Roughly, one has
\begin{equation}\aligned\label{adH4a}
d(\alpha;y)\leq2
.
\endaligned\end{equation}

Therefore, by \eqref{adH2}, \eqref{adH3} and \eqref{adH4a},
\begin{equation}\aligned\label{adH5}
&\sum_{n,m}(n^2-\frac{(m+\frac{n}{2})^2}{y^2})^2 e^{-\pi\alpha(yn^2+\frac{(m+\frac{n}{2})^2}{y})}\\
\leq&
4(1-\frac{1}{4y^2})^2(y+\frac{1}{4y})^2 e^{-\pi\alpha(y+\frac{1}{4y})}
+\frac{2}{y^6}e^{-\frac{\pi\alpha}{y}}
+2\cdot 16^2 e^{-4\pi\alpha y}\\
&\;\;\;\;+2\cdot16^2 y^2 e^{-\pi\alpha(4y+\frac{1}{y})}.
\endaligned\end{equation}
Then the desired result follows.

\end{proof}

\section{Proof of Theorems \ref{Th1}-\ref{Th3}}

{\bf Proof of Theorem \ref{Th1}.}
\vskip0.05in
Case 1: $b=\frac{1}{2\pi}$. This follows from Theorems \ref{Th31} and \ref{Th41}.

Case 2: $b<\frac{1}{2\pi}$. It is proved by Lemma \ref{Lemma1} and Case 1.

Case 3: $b>\frac{1}{2\pi}$. It follows by Lemma \ref{Lemma32}. Indeed, by Lemma \ref{Lemma32}, one has
 \begin{equation}\aligned\nonumber
\mathcal{W}_b(\alpha;z)&=\alpha^{-\frac{3}{2}}\sqrt y\cdot\Big(\frac{1}{2\pi}-b+o(1)
\Big)\\
&\mapsto -\infty,\;\;\;\hbox{if}\;\;\;b>\frac{1}{2\pi},\;\;\hbox{as}\;\; y\rightarrow+\infty
\endaligned\end{equation}
proves the nonexistence result.

\vskip0.05in

{\bf Proof of Theorem \ref{Th3}.}
\vskip0.05in
Case 1: $b=\sqrt{a}$. By simple observation, one has the connection between the functional $\theta(\alpha;z)-\sqrt a\theta(a\alpha;z)$
and $\mathcal{W}_{\frac{1}{2\pi}}(x\alpha;z)$. Indeed,
applying the fundamental Theorem of calculus on a parameter $t$, we have

 \begin{equation}\aligned\label{W1}
\theta(\alpha;z)-\sqrt a\theta(a\alpha;z)=-\int_1^a \partial_t(\sqrt t\theta(t\alpha;z))dt
=\pi\int_1^a\mathcal{W}_{\frac{1}{2\pi}}(t\alpha;z)dt.
\endaligned\end{equation}
See $\mathcal{W}_{\frac{1}{2\pi}}(t\alpha;z)$ in \eqref{Wbbb}.
The proof then follows by Theorem \ref{Th1}(or Theorem \ref{Th31}) and \eqref{W1}.

Case 2: $b<\sqrt{a}$.

\begin{equation}\aligned\label{W2}
\Big(\theta(\alpha;z)-b\theta(a\alpha;z)\Big)
=\Big(\theta(\alpha;z)-\sqrt a\theta(a\alpha;z)\Big)
+(\sqrt a-b)\theta(a\alpha;z).
\endaligned\end{equation}
Then it follows by \eqref{W2}, Case 1 and  \begin{equation}\aligned\label{GH1}
\min_{z\in\mathbb{H}}\theta(\alpha;z)\;\;\hbox{is achieved at}\;\; e^{i\frac{\pi}{3}}
\endaligned\end{equation}
by \cite{Mon1988}.

Case 3: $b>\sqrt{a}$.
By Lemma \ref{Lemma34}, for $\forall\alpha>0$,
\begin{equation}\aligned\nonumber
\theta (\alpha;z)-b\theta(a\alpha;z)
&=\sqrt{\frac{y}{a\alpha}}\cdot\Big(
\sqrt a-b+o(1)
\Big)\\
&\mapsto -\infty,\;\;\;\hbox{if}\;\;\;b>\sqrt a,\;\;\hbox{as}\;\; y\rightarrow+\infty
\endaligned\end{equation}
which proves the nonexistence result.

\bigskip
\noindent
{\bf Acknowledgements.}
 The research of S. Luo is partially supported by NSFC(Nos. 12261045, 12001253) and double thousands plan of Jiangxi(jxsq2019101048). The research of J. Wei is partially supported by NSERC of Canada.

{\bf Statements and Declarations: there is no conflict of interest.}

{\bf Data availability: the manuscript has no associated data.}

\end{document}